\documentclass[12pt]{article}
\usepackage[latin1]{inputenc}
\usepackage{amssymb,amsmath,amsfonts}
\usepackage{xcolor}
\newtheorem{theorem}{Theorem}
\newtheorem{proposition}[theorem]{Proposition}
\newtheorem{lemma}[theorem]{Lemma}

\newtheorem{corollary}[theorem]{Corollary}

\newtheorem{remark}[theorem]{Remark}

\newtheorem{example}[theorem]{Example}

\newcommand{\Ce}{\mathcal{C}}
\newcommand{\R}{\mathbb{R}}

\newcommand{\Q}{\mathbb{Q}}
\newcommand{\Sf}{\mathbb{S}}
\newcommand{\C}{\mathbb{C}}

\newcommand{\spa}{\mbox{span}}

\newcommand{\grad}{\mbox{grad}}

\newcommand{\po}{{\hspace*{-1ex}}{\bf .  }}
\newcommand{\ii}{isometric immersion }
\newcommand{\iis}{isometric immersions }

\newcommand{\fa}{\text{for all \,\ }}
\newcommand{\h}{\mathcal{H}}
\newcommand{\1}{\phi(e_r)}
\newcommand{\2}{\omega(e_r)}

\def\t{{\theta}}
\def\gh{{\hat{g}}}

\def\D{{\Delta}}

\def\<{\langle}
\def\n{\nabla}

\def\>{\rangle}
\def\a{\alpha}
\def\v{\varphi}

\def\bea{\begin{eqnarray*} }
\def\eea{\end{eqnarray*} }
\def\be{\begin{equation}}
\def\ee{\end{equation}}

\def\nap{\nabla^\perp}

\def\proof{\noindent{\it Proof: }}
\def\qed{\ifhmode\unskip\nobreak\fi\ifmmode\ifinner
\else\hskip5 pt \fi\fi\hbox{\hskip5 pt \vrule width4 pt
height6 pt  depth1.5 pt \hskip 1pt }}
\begin{document}
\title{On the mean curvature of submanifolds with nullity}
\author{A.E.\ Kanellopoulou and Th.\ Vlachos*}
\date{}
\maketitle
\renewcommand{\thefootnote}{\arabic{footnote}} 

\renewcommand{\thefootnote}{\fnsymbol{footnote}} 
\renewcommand{\thefootnote}{\arabic{footnote}}

\renewcommand{\thefootnote}{\fnsymbol{footnote}} 
\renewcommand{\thefootnote}{\arabic{footnote}}

\renewcommand{\thefootnote}{\fnsymbol{footnote}} 
\footnotetext{The first named author would like to acknowledge financial support by the General Secretariat for Research and Technology (GSRT) and the Hellenic Foundation for Research and Innovation (HFRI), Grant No: 133, Action: Support for Postdoctoral Researchers.}  
\renewcommand{\thefootnote}{\arabic{footnote}}  

\begin{abstract}
 In this paper, we investigate geometric conditions for isometric immersions  
with positive index of relative nullity  to be cylinders.  There is an abundance of noncylindrical $n$-dimensional minimal submanifolds with index of relative nullity $n-2$,  fully described by Dajczer and Florit  \cite{DF2} in terms of a  certain class of elliptic surfaces. Opposed to this, we prove that nonminimal $n$-dimensional submanifolds in space  forms of any codimension are locally cylinders provided that they carry a totally geodesic distribution of rank $n-2\geq2,$ which is contained in the relative nullity distribution, such that   the length of the mean curvature vector field is constant along each leaf.
The case of dimension $n=3$ turns out to be special. We show that  there exist elliptic three-dimensional submanifolds in spheres satisfying the above properties. In fact,  we provide a parametrization of three-dimensional submanifolds as unit tangent bundles of minimal surfaces in the Euclidean space whose first curvature ellipse is nowhere a circle and its second one is everywhere a circle. Moreover, we provide several applications to submanifolds whose mean curvature vector field has constant length, a much weaker condition than being parallel. 
\end{abstract}

 \section{Introduction}

A fundamental concept in the theory of submanifolds is the index of relative nullity introduced by Chern and Kuiper \cite{chrn}. At a point $x\in M^n$ the \textit{index of relative nullity} $\nu(x)$ of an isometric immersion $f\colon M^n\to\Q_c^m$  is the dimension of the \textit{relative nullity} tangent subspace $\D_f(x)$ of $f$ at $x$, that is, 
the kernel of the second fundamental form $\alpha^f$ at that point. Here, $\mathbb{Q}_c^m$ is the simply connected space form with curvature $c$, that is,  the Euclidean space $\R^m,$ the sphere $\mathbb{S}^m$ or the hyperbolic space $\mathbb{H}^m$, according to whether $c=0,c=1$ or $c=-1,$ respectively.
The kernels form an integrable distribution along any open subset where 
the index is constant and the images under $f$  of the leaves of the foliation  
are totally geodesic submanifolds in the ambient space. 

Cylinders are the simplest examples of submanifolds with positive index of relative nullity. An isometric immersion $f\colon M^n\to \R^m$ is said to be a $k$-\textit{cylinder} if the manifold $M^n$ splits as a Riemannian product 
$M^n=M^{n-k}\times\R^k$ and there is an isometric immersion $g\colon M^{n-k}\to\R^{m-k}$ such that
$f=g\times{\rm{id}}_{\R^k}$.
A natural problem in submanifold theory is to find geometric conditions 
for an isometric immersion with index of relative nullity $\nu\geq k>0$ at 
any point to be a $k$-cylinder. 

A fundamental result asserting that an isometric immersion 
$f\colon M^n\to\R^m$ of a complete Riemannian manifold with positive index of relative nullity  must be a 
cylinder is Hartman's theorem \cite{har} that requires the Ricci 
curvature of $M^n$ to be nonnegative.
Even for hypersurfaces, the same conclusion does not hold  if instead we 
assume that the Ricci curvature is nonpositive. 
Notice that the latter is always the case if $f$ is a minimal 
immersion. Counterexamples easy to construct are the complete irreducible 
ruled hypersurfaces of any dimension discussed in \cite[p.\ 409]{dg}. 

The cylindricity of minimal submanifolds was studied in \cite{D, hsv} under global assumptions.
These results are truly global in nature since there are plenty of (noncomplete) 
examples of minimal submanifolds of any dimension $n$ with constant index $\nu=n-2$ 
that are not part of a cylinder on any open subset. 
They can be all locally parametrically  described in terms of a certain 
class of elliptic surfaces (see \cite[Th. 22]{DF2}). 
Some of the many papers containing  characterizations of submanifolds as 
cylinders without the requirement of minimality are \cite{dg0,gf,har,ma}. 

In this paper, we deal with nonminimal $n$-dimensional submanifolds  of arbitrary codimension and index of relative nullity $\nu\geq n-2$ at any point. Our aim is to provide geometric conditions, in terms of the mean curvature, for an isometric immersion to be a cylinder. The choice of the geometric condition is inspired by the observation that cylinders are endowed with a totally geodesic distribution contained in the relative nullity distribution, such that the mean curvature is constant along each leaf.
 Throughout the paper,  the \textit{mean curvature} of an isometric immersion $f$ is defined as the length $H=\| \h \|$ of the \textit{mean curvature vector field} given by $\h=\mathrm{trace}(\alpha^f)/n.$ 
 
The following result provides a characterization of cylinders of dimension $n\geq4.$ 

\begin{theorem}\po\label{main}
Let $f\colon M^n \to \mathbb{Q}^{n+p}_c, n\geq4,$ be an isometric immersion such that $M^n$ carries a totally geodesic distribution $D$ of rank $n-2$ satisfying $D(x)\subseteq \D_f(x)$ for any $x\in M^n$. If the mean curvature of $f$ is constant along each leaf of $D$, then either $f$ is minimal or $c=0$ and $f$ is locally a $(n-2)$-cylinder over a surface on the open subset where the mean curvature is positive. Moreover, the submanifold is globally a cylinder if  the leaves of $D$ are complete. 
\end{theorem}

It is interesting that the above theorem fails for substantial three-dimensional submanifolds of codimension $p\geq2$. Being substantial means that the codimension cannot be reduced. We show that besides cylinders, there exist elliptic three-dimensional submanifolds in spheres satisfying the properties assumed in Theorem \ref{main}.  Thus the submanifolds being three-dimensional is special. The notion of elliptic submanifolds was introduced in \cite{DF2}.
In fact,  the following result allows a parametrization of them in terms of minimal surfaces in the Euclidean space, the so-called \textit{bipolar parametrization}, using the following construction.

Let $g\colon L^2 \to \R^{n+1}, n\geq5,$ be a minimal surface. The map 
$\Phi_g\colon T^1L\to \mathbb{S}^n$ defined on the unit tangent bundle of $L^2$ and given by 
\be \label{start}
\Phi_g(x,w)=g_{*_x}w
\ee
parametrizes (outside singular points) an immersion with index of relative nullity at least one at any point. 

\begin{theorem}\po\label{main1}
Let $f\colon M^3 \to \mathbb{Q}^{3+p}_c$ be an isometric immersion such that $M^3$ carries a totally geodesic distribution $D$ of rank one satisfying $D(x)\subseteq \D_f(x)$ for any $x\in M^3$. If the mean curvature of $f$ is constant along each integral curve of $D$, then one of the following holds:

(i) The immersion $f$ is minimal.

(ii)  $c=0$ and  $f$ is locally a cylinder over a surface. 

(iii)  $c=1$ and the immersion $f$ is elliptic and locally parametrized by \eqref{start}, where $g\colon L^2 \to \R^{n+1}, n\geq5,$ is a minimal surface whose first curvature ellipse is nowhere a circle and the second curvature ellipse is everywhere a circle. 
\end{theorem}

Minimal surfaces satisfying the conditions in part (iii) of the above theorem can be constructed using the Weierstrass representation  by choosing appropriately the holomorphic data. It is worth noticing that minimal surfaces in the Euclidean space that satisfy the Ricci condition, or equivalently are locally isometric to a minimal surface in $\R^3,$ fulfill these conditions (see Section 6 for details). These surfaces were classified by Lawson \cite{Law}.

The above results allow us to provide applications to submanifolds with constant mean curvature and not necessarily constant positive index of relative nullity.

Having constant mean curvature is a much weaker restriction on the mean curvature vector field than being parallel in the normal bundle. One can check that three-dimensional elliptic submanifolds described in Theorem \ref{main1} do not have parallel mean curvature vector field along the totally geodesic distribution. Combining this with Theorem \ref{main}, it follows that a submanifold is locally a cylinder provided that it carries a totally geodesic distribution of rank $n-2\geq1$ that is  contained in the relative nullity distribution, along  which the mean curvature vector field is parallel in the normal connection.

Opposed to the fact that there is an abundance of noncylindrical $n$-dimensional minimal submanifolds with index of relative nullity $n-2$ (see \cite{DF2}),  we prove the following result for submanifolds with constant positive mean curvature.

\begin{theorem}\po\label{main2}
Let $f\colon M^n \to \mathbb{Q}^{n+p}_c,n\geq3,$ be a nonminimal isometric immersion with index of relative nullity $\nu\geq n-2$ at any point. If  the mean curvature of $f$ is constant and either $n\geq4$ or $n=3$ and $p=1$, then $c=0$. Moreover, there exists an open dense subset $V\subseteq M^n$ such that every point has a neighborhood $U\subseteq V$ so that $f(U)$ is an open subset of the image of a cylinder either over a surface in $\R^{p+2}$, or over a curve in $\R^{p+1}$ with  constant first Frenet curvature.
\end{theorem}

The following is an immediate consequence of the above result due to real analyticity of hypersurfaces with constant mean curvature.  
\begin{corollary}\po\label{main2'}
Let $f\colon M^n \to \mathbb{Q}^{n+1}_c, n\geq3,$ be a nonminimal isometric immersion with index of relative nullity $\nu\geq n-2.$ If  the mean curvature of $f$ is constant, then $c=0$ and $f(M)$ is an open subset of the image of a cylinder over a surface in $\R^3$ of constant mean curvature.
\end{corollary}

The next result extends Corollary 1 in \cite{cfgms} for hypersurfaces in every space form without any global assumption.

\begin{corollary}\po\label{T5}
Let $f\colon M^n \to \mathbb{Q}^{n+1}_c, n\geq3,$ be an isometric immersion with constant mean curvature. If $M^n$ has sectional curvature $K\leq c$, then either $f$ is minimal or $c=0$ and $f(M)$ is an open subset of the image of a cylinder over a surface in $\R^3$ of constant mean curvature.
In the latter case, $f$ is a cylinder over a circle provided that $M^n$ is complete. 
\end{corollary}

The following rigidity result that was proved in \cite{dg0} for $c=0$ is another consequence of our main results.

\begin{corollary}\po\label{C6}
Any nonminimal isometric immersion $f\colon M^n \to \mathbb{Q}^{n+1}_c,n\geq3,$ with constant mean curvature is  rigid, unless $c=0$ and $f(M)$ is an open subset of the image of a cylinder over a surface in $\R^3$ of constant mean curvature.
\end{corollary}

Our next result extends to any dimension a well-known theorem for constant mean curvature surfaces due to Klotz and Osserman \cite{ko} (see \cite{Luis} for another extension). 

\begin{theorem}\po\label{T7}
Let $f\colon M^n \to \mathbb{Q}_c^{n+1},n\geq3,$ be an isometric immersion with constant mean curvature,  where $c=0$ or $c=1$. If $M^n$ is complete and its extrinsic curvature does not change sign, then either $f$ is minimal or  totally umbilical or a cylinder over a sphere of dimension $1\leq k < n.$   
\end{theorem}

For submanifolds with constant mean curvature of codimension two, we prove the following.

\begin{theorem}\po\label{T8}
Let $f\colon M^n \to \R^{n+2}, n\geq3,$ be a nonminimal isometric immersion with constant mean curvature. If  the sectional curvature of $M^n$ is nonpositive, then there exists an open dense subset $V\subseteq M^n$ such that  every point has a neighborhood $U\subseteq V$ where one of the following holds:

(i) The neighborhood $U$ splits as a Riemannian product  $U=M^2\times W^{n-2}$ such that  $f|_U=g\times j$ is a product, where  $g\colon M^2 \to\R^4$ is a surface with constant mean curvature and $j\colon W^{n-2}\to \R^{n-2}$ is the inclusion. 

 (ii) The immersion on $U$ is a composition  $f|_U=h\circ F,$ where $h=\gamma\times id_{\R^{n-1}} \colon  \R\times\R^n\to\R^{n+2}$ is cylinder over a unit speed plane curve $\gamma(s)$ with curvature $k(s)$ and $F\colon M^n \to\R^{n+1}$ is a hypersurface.  Moreover, the mean curvature $H_F$ of $F$ is given by 
 $$
 H_F^2=H_f^2-\frac{1}{n^2}k^2\circ F_a\left(1-\<\xi, a\>^2\right)^2,
 $$
 where $F_a=\<F,a\>$ and $\<\xi, a\>$ are the height functions of $F$ and its Gauss map $\xi$ relative to the unit vector $a=\partial/\partial s,$  respectively.

 (iii) The neighborhood $U$ splits as a Riemannian product  $U=M^2_1\times M_2^2\times W^{n-4}$ such that $f|_U=g_1\times g_2\times j$ is a product, where $g_i\colon M^2_i \to\R^3, i=1,2,$ are surfaces with constant mean curvature and  $j\colon W^{n-4}\to \R^{n-4}$ is the inclusion.  
\end{theorem}

For constant sectional curvature submanifolds with constant mean curvature of codimension two, we prove the following theorem that extends results in \cite{dtt, enomoto1}.

\begin{theorem}\po\label{T11}
Let $f\colon M_{\tilde{c}}^n\to\mathbb{Q}_c^{n+2},n\geq3,$ be an isometric immersion of a Riemannian manifold of constant sectional curvature $\tilde c.$ If the mean curvature of $f$ is constant and either $n\geq4$ or $n=3$ and $c=\tilde{c},$ then one of the following holds:
 
(i) $f$ is totally geodesic or totally umbilical.

(ii) $\tilde{c}=c=0$ and $f=g\times j,$ where $g\colon M^2\to\R^4$ is a flat surface with constant mean curvature and  $j\colon W\to\R^{n-2}$ is an inclusion.

(iii) $\tilde{c}=0, c=-1$ and $f$ is a composition $f=i\circ F$, where $i\colon \R^{n+1}\to \mathbb H^{n+2}$ is the inclusion as a horosphere and  $F\colon M^n_{\tilde{c}} \to \R^{n+1}$ is cylinder over a circle.
\end{theorem}

Cylinder theorems for complete minimal K\"ahler submanifolds were proved in \cite{dr,fur}. For K\"ahler submanifolds with constant mean curvature, we prove the following results.

\begin{theorem}\po\label{T9}
Let $f\colon M^n \to \R^{n+1},n\geq4,$ be an isometric immersion with constant mean curvature. If $M^n$ is K\"ahler, then either $f$ is minimal or $f(M)$ is an open subset of the image of a cylinder over a surface in $\R^3$ with constant mean curvature.
\end{theorem}

\begin{theorem}\po\label{T10}
Let $f\colon M^n \to \R^{n+2}, n\geq4,$ be a nonminimal isometric immersion of a K\"ahler manifold $M^n$  with constant mean curvature. If  the Ricci curvature or the holomorphic curvature of $M^n$ is nonnegative, then there exists an open dense subset $V\subseteq M^n$ such that every point has a neighborhood $U\subseteq V$ where $f|_U$ is as in Theorem \ref{T8}.
\end{theorem}

The paper is organized as follows: In Section 2, we recall well-known results about the relative nullity distribution, totally geodesic distributions that are contained in the relative nullity distribution, as well as results about their splitting tensor.
In Section 3, we fix the notation, give some preliminaries and prove auxiliary results that will be used in the proofs of our main theorems.  Section 4 is devoted to the proof of Theorem \ref{main}. 
In Section 5, we recall the notion of elliptic submanifolds, as well as the associated notion of higher curvature ellipses. 
We also discuss the polar and bipolar surfaces of elliptic submanifolds.  
In Section 6, we study the case of three-dimensional submanifolds. We provide a parametrization for these submanifolds in terms of certain elliptic surfaces, the so-called \textit{polar parametrization} (see Theorem \ref{main0}). Based on this, we give the proof of Theorem \ref{main1}. We conclude this section by showing that minimal surfaces in the Euclidean space that  are locally isometric to a minimal surface in $\R^3$ satisfy the conditions in part (iii) of Theorem \ref{main1}.  In Section 7, we prove Theorem \ref{main2} and the applications of our main results on submanifolds with constant mean curvature. In addition, we provide examples of submanifolds as in part (ii) of Theorems \ref{T8} and \ref{T11}.

\section{The relative nullity distribution}

In this section, we recall some basic facts from the theory of \iis
that will be used throughout the paper.

Let $M^n,n\geq3,$ be a Riemannian manifold and let $f\colon M^n\to\Q^m_c$ be an isometric immersion into a space form $\Q^m_c$. The \emph{relative nullity} 
subspace $\Delta_f(x)$ of $f$ at $x\in M^n$ is the  kernel of its second fundamental 
form $\alpha^f\colon TM\times TM\to N_fM$ with values in the normal bundle, 
that is,
$$
\Delta_f(x)=\left\{X\in T_xM:\alpha^f(X,Y)=0\;\;\text{for all}\;\;Y\in T_xM\right\}.
$$
The dimension $\nu(x)$ of $\Delta_f(x)$ is called the \emph{index of relative nullity}  
of $f$ at $x\in M^n$. 

A smooth distribution $D\subset TM$ on $M^n$ is \emph{totally geodesic} if $\n_TS\in \Gamma(D)$ whenever $T,S\in \Gamma(D)$.
Let $D$ be a smooth distribution on $M^n$ and $D^{\perp}$ denote the distribution on $M^n$ that assigns to each $x\in M^n$ the orthogonal complement of $D(x)$ in $T_xM$. We write $X=X^v+X^h$
according to the orthogonal splitting $TM=D\oplus D^{\perp}$
and denote ${\nabla}^{h}_XY = (\nabla_X Y)^h$ for all $X,Y\in TM$, where $\n$ is the Levi-Civit\'a connection on $M^n.$
The \textit{splitting tensor} $C\colon D\times D^{\perp}\to D^{\perp}$ 
is given by
$$
C(T,X)=-{\nabla}^{h}_XT
$$
for any $T\in D$ and $X\in D^{\perp}$. 

When $D$ is a totally geodesic distribution such that  $D(x)\subseteq \D_f(x)$ for all $x\in M^n,$ the following differential equation for the tensor $\Ce_T=C(T,\cdot)$
is well-known to hold (cf. \cite{dg} or  \cite{da}):
\be\label{C1}
\nabla^{h}_S \Ce_T=\Ce_T\circ \Ce_S+\Ce_{\n_ST} +c\<S,T\>I,
\ee
where $I$ is the identity endomorphism of $D^{\perp}$. Here $\nabla^{h}_S \Ce_T \in \Gamma ({\rm{End}} (D^\perp))$ is defined by
$$
(\nabla^{h}_S \Ce_T)X=\nabla^{h}_S \Ce_TX-\Ce_T\nabla^{h}_SX
$$
for all $T,S\in D$ and $X\in D^{\perp}$. The Codazzi equation gives
\be\label{cod}
\n_TA_{\xi}=A_{\xi}\circ\Ce_T+ A_{\n^{\perp}_T\xi}
\ee
for any $T\in D,$ where the shape operator $A_{\xi}$ with respect to the normal direction $\xi$ is restricted to 
$D^{\perp}$ and $\n^{\perp}$ stands for the normal connection of $f$.
In particular, the endomorphism $A_\xi\circ\Ce_T$ of $D^\perp$ is symmetric, that is,
\be\label{C3}
A_{\xi}\circ\Ce_T=\Ce_T^t\circ A_{\xi}. 
\ee

For later use, we recall the following known results. 

\begin{proposition} (\cite[Prop. 7.4]{da}) \po\label{cylinder}  Let $f\colon M^n\to\mathbb{Q}^m_c$ be 
an isometric immersion such that $M^n$ carries a smooth totally geodesic distribution $D$ of rank $0<k< n$ satisfying $D(x)\subseteq \D_f(x)$ for all $x\in M^n.$ If the splitting tensor $C$
vanishes, then $c=0$ and $f$ is locally a $k$-cylinder.
\end{proposition}

\begin{proposition} (\cite[Prop. 1.18]{da}) \po\label{nu0}
For an isometric immersion $f\colon M^n\to\Q_c^m$, the following assertions hold:

(i) The index of relative nullity $\nu$ is upper semicontinuous. In particular, the subset 
$$
M_0=\{x\in M^n:\nu(x)=\nu_0\},
$$
where $\nu$ attains its minimum value $\nu_0$ is open.

(ii) The relative nullity distribution $x\mapsto\Delta_f(x)$ is smooth on any subset of $M^n$ where $\nu$ is constant.

(iii) If $U\subseteq M^n$ is an open subset where $\nu$ is constant, then $\D_f$ is a totally geodesic (hence integrable) distribution on $U$ and the restriction of $f$  to each leaf is totally geodesic.  
\end{proposition}

 \section{Auxiliary results}
The aim of  this section is to prove several lemmas that will be used in the proofs of our main results. 
 
Throughout this section, we assume that $f\colon M^n \to \mathbb{Q}^{n+p}_c,n\geq3,$ is a nonminimal isometric immersion such that $M^n$ carries a smooth totally  geodesic distribution $D$ of rank $n-2$ satisfying $D(x)\subseteq \D_f(x)$ for any $x\in M^n$. We also assume that the mean curvature of $f$ is constant along each leaf of $D$.

Hereafter we work on the open subset where the mean curvature is positive and choose a local orthonormal frame
$\xi_{n+1}, \dots, \xi_{n+p}$
in the normal bundle $N_fM$, such that $\xi_{n+1}$ is collinear to the mean curvature vector field. We also choose a local orthonormal frame $e_1, \dots, e_n$ in the tangent bundle $TM$ such that $e_1, e_2$ span $D^\perp$  and diagonalize $A_{\xi_{n+1}}|_{D^\perp}$, where $A_{\xi_{n+1}}$ denotes the shape operator of $f$ with respect to 
$\xi_{n+1}$. Then, we have $A_{\xi_{n+1}}e_i=k_ie_i, \,\ i=1,2,
$ and  consequently the mean curvature is given by $nH=k_1 +k_2,$ where $k_1,k_2$ are the principal curvatures.

Since the mean curvature is positive, at least one of the principal curvatures $k_1$ and $k_2$ has to be different from zero. In the sequel, we assume without loss of generality, that $k_1\neq 0$ and define the function
 $$
\rho=-\frac{k_2}{k_1}.
$$
On the open subset where the mean curvature is positive we have
\begin{equation}\label{k1k2}
k_1=-\frac{nH}{\rho-1}\,\  \,\ \text{and}\,\ \,\ k_2=\frac{n\rho H}{\rho-1}.
\end{equation}
The above mentioned notation is used throughout the paper.

The following lemma gives the form of the splitting tensor.

\begin{lemma}\po\label{char}
On the open subset where the mean curvature is positive, the splitting tensor is given by
$$
\Ce_T=\psi_1(T)L_1 + \psi_2(T)L_2
$$
for any $T\in\Gamma(D),$ where $\psi_1, \psi_2$ are 1-forms dual to the vector fields $\n_{e_2}e_2, \n_{e_1}e_2,$ respectively, and $L_1$, $L_2$ $\in\Gamma({\rm{End}}(D^\perp))$ are defined by $L_1e_1=\rho e_1=-L_2e_2$ and $L_1e_2=e_2 =L_2e_1.$
Moreover, the following holds:
\begin{align}
& T(k_1)=\rho k_1\psi_1(T)+\sum_{\alpha=n+2}^{n+p}\<\nap_{T}\xi_{n+1},\xi_{\alpha}\>\<A_{\xi_\alpha}e_1,e_1\>, \label{1}\\
& T(k_2)=k_2\psi_1(T)-\sum_{\alpha=n+2}^{n+p}\<\nap_{T}\xi_{n+1},\xi_{\alpha}\>\<A_{\xi_\alpha}e_1,e_1\>, \label{2}\\
& (k_1-k_2)\omega(T)=k_2\psi_2(T) + \sum_{\alpha={n+2}}^{n+p}\<\nap_{T}\xi_{n+1},\xi_{\alpha}\>\<A_{\xi_\alpha}e_1,e_2\>, \label{3} \\
& (k_1-k_2)\omega(T)=-\rho k_1\psi_2(T)+ \sum_{\alpha={n+2}}^{n+p}\<\nap_{T}\xi_{n+1},\xi_{\alpha}\>\<A_{\xi_\alpha}e_1,e_2\> \label{4}
\end{align}
for any $T\in\Gamma(D)$, where  $\omega$ denotes the connection form given by $\omega=\<\n e_1,e_2\>$. 
\end{lemma}

\proof 
From the Codazzi equation  we have
$$
\big(\n_{T}A_{\xi_{n+1}}\big)e_i - \big(\n_{e_i}A_{\xi_{n+1}}\big)T=A_{\nap_T{\xi_{n+1}}}e_i - A_{\nap_{e_i}{\xi_{n+1}}}T
$$
for any $T\in\Gamma(D)$ and $i=1,2.$ The above is equivalent to the following:
\begin{align*}
& T(k_1)=k_1\<\n_{e_1}e_1,T\>+\sum_{\alpha=n+2}^{n+p}\<\nap_{T}\xi_{n+1},\xi_{\alpha}\>\<A_{\xi_\alpha}e_1,e_1\>, \\ 
& T(k_2)=k_2\<\n_{e_2}e_2,T\>-\sum_{\alpha=n+2}^{n+p}\<\nap_{T}\xi_{n+1},\xi_{\alpha}\>\<A_{\xi_\alpha}e_1,e_1\>, \\
& (k_1-k_2)\omega(T)=k_2\<\n_{e_1}e_2,T\> + \sum_{\alpha={n+2}}^{n+p}\<\nap_{T}\xi_{n+1},\xi_{\alpha}\>\<A_{\xi_\alpha}e_1,e_2\>,  \\
& (k_1-k_2)\omega(T)=k_1\<\n_{e_2}e_1,T\>+ \sum_{\alpha={n+2}}^{n+p}\<\nap_{T}\xi_{n+1},\xi_{\alpha}\>\<A_{\xi_\alpha}e_1,e_2\>.
\end{align*}

Using the assumption that the mean curvature is constant along each leaf of the distribution $D$, the first two equations imply
$$
\<\n_{e_1}e_1,T\>=\rho\<\n_{e_2}e_2,T\>
$$
for any $T\in\Gamma(D).$
Additionally, the last two equations yield
$$
\<\n_{e_2}e_1,T\>=-\rho\<\n_{e_1}e_2,T\>.
$$
Now the structure of the splitting tensor and \eqref{1}-\eqref{4} follow easily from the above. \qed

\begin{lemma}\po\label{eq}
Let $e_r, r\geq3,$ be an orthonormal frame of the distribution  $D.$ Then the functions $u_r:=\psi_1(e_r)$ and $v_r:=\psi_2(e_r)$ satisfy 
\begin{equation}
2\rho(u_ru_s + v_rv_s)-c\delta_{rs} = \frac{\rho-1}{nH}\sum_{\alpha=n+2}^{n+p}\<\nap_{e_r}\xi_{n+1},\xi_{\alpha}\>\big(u_s\<A_{\xi_\alpha}e_1,e_1\> -v_s\<A_{\xi_\alpha}e_1,e_2\>  \big)\label{second}
\end{equation}
for all $ r, s\geq3$, where $\delta_{rs}$ is the Kronecker delta.
\end{lemma}

\proof Using Lemma \ref{char}, we have
\begin{align}
(\n^h_{e_r}\Ce_{e_s})=e_r(u_s)L_1+e_r(v_s)L_2+u_s\n^h_{e_r}L_1+v_s\n^h_{e_r}L_2\label{dc}
\end{align}
 for any $r,s\geq3$. A direct computation yields
\begin{align}
&(\n^h_{e_r}L_1)e_1=-(\n^h_{e_r}L_2)e_2=e_r(\rho)e_1+(\rho-1)\2e_2,\label{dl11}\\
&(\n^h_{e_r}L_1)e_2=(\n^h_{e_r}L_2)e_1=(\rho-1)\2e_1.\label{dl12}
\end{align} 
Then \eqref{1} and \eqref{2} imply that
\be\label{ero}
e_r(\rho)=-\rho(\rho-1)u_r + \frac{(\rho-1)^2}{nH}\sum_{\alpha=n+2}^{n+p}\<\nap_{e_r}\xi_{n+1},\xi_{\alpha}\>\<A_{\xi_\alpha}e_1,e_1\>.
\ee
From (\ref{C1}) we know that the splitting tensor satisfies
\begin{align}
&(\n_{e_r}^h\Ce_{e_s})e_i=\Ce_{e_s}\circ\Ce_{e_r}e_i + \Ce_{\n_{e_r}e_s}e_i +c\delta_{rs}e_i \label{crs1}
\end{align}
for any $r,s\geq3$ and $i=1,2$.

Let $\omega_{rs}$  be the connection form given by 
$\omega_{rs}=\<\n e_r,e_s\> $ for all $r,s\geq3$. Using \eqref{dc}-\eqref{ero}, we find that \eqref{crs1} for $i=1$ is equivalent to
\begin{align}  \label{rv1}
\rho e_r(u_s)=&\rho(2\rho-1)u_ru_s-\rho v_rv_s -(\rho-1)v_s\omega(e_r)\qquad \nonumber\\ 
&-u_s\frac{(\rho-1)^2}{nH} \sum_{\alpha=n+2}^{n+p}\<\nap_{e_r}\xi_{n+1},\xi_{\alpha}\>\<A_{\xi_\alpha}e_1,e_1\>
+ \rho\sum_{t\geq3}^n\omega_{st}(e_r)u_t+c\delta_{rs}
\end{align}
and
\begin{align}
e_r(v_s)=\rho u_r v_s +u_s v_r-(\rho-1)u_s\2+\sum_{t\geq3}^n\omega_{st}(e_r)v_t \qquad \qquad \qquad \qquad\label{ru}
\end{align}
for all $r,s\geq3$. Moreover, \eqref{crs1} for $i=2$ implies that
\begin{align}
e_r(u_s)=&u_ru_s-\rho v_rv_s+(\rho-1)v_s\2+\sum_{t\geq3}^n\omega_{st}(e_r)u_t+c\delta_{rs} \qquad \qquad \qquad  \label{rv2}
\end{align}
for all $r,s\geq3$. 

Combining \eqref{rv1} and \eqref{ru},  we obtain 
$$
2\rho u_ru_s + \rho v_rv_s - c\delta_{rs}- v_s(\rho+1)\2=u_s\frac{\rho-1}{nH} \sum_{\alpha=n+2}^{n+p}\<\nap_{e_r}\xi_{n+1},\xi_{\alpha}\> \<A_{\xi_\alpha}e_1,e_1\>. 
$$
Using \eqref{k1k2}, it is easily seen that \eqref{4} is written as
\be\label{omega}
(\rho+1)\2=-\rho  v_r-\frac{\rho-1}{nH} \sum_{\alpha=n+2}^{n+p}\<\nap_{e_r}\xi_{n+1},\xi_{\alpha}\>\<A_{\xi_\alpha}e_1,e_2\>,
\ee
and now \eqref{second} follows directly from the above two equations. \qed
\vspace{2ex}

We recall that the \textit{first normal space} $N_1^f(x)$ of the immersion $f$ at a point $x\in M^n$ is the subspace of its normal space $N_fM(x)$ spanned by the image of its second fundamental form $\alpha^f$ at $x$, that is,
$$
N_1^f(x)={\rm{span}}\left\{\alpha^f(X,Y): X,Y \in T_xM\right\}.
$$
The rank condition and the symmetry of the second fundamental form imply that $\dim N_1^f(x)\leq3$ for all $x\in M^n$. 

Consider the open subset 
$$
M_3=\left\{x\in M^n:\dim N_1^f(x)=3\right\}.
$$

\begin{lemma}\po\label{M3}
The splitting tensor vanishes on the open subset $M_3^*:=M_3\smallsetminus\{x\in M^n: H(x)=0\}$.
\end{lemma}

\proof
On the subset $M_3^*,$ we consider the orthogonal splitting $
N_1^f=\hat{N}^f_1 \oplus {\rm{span}}\{\h\}.$
Choose the local frame such that $\xi_{n+1}$ is collinear to the mean curvature vector field $\h$, and  $ \xi_{n+2},\xi_{n+3}$ span the plane bundle $\hat{N}^f_1$. Then, we have
 $$
 \text{trace}A_{\xi_{n+2}}|_{D^\perp}=0=\text{trace}A_{\xi_{n+3}}|_{D^\perp}.
 $$
 Hence, we obtain
$$
A_{\xi_{n+2}}|_{D^\perp}\circ J=J^t\circ A_{\xi_{n+2}}|_{D^\perp} \,\ \text{and} \,\ A_{\xi_{n+3}}|_{D^\perp}\circ J=J^t\circ A_{\xi_{n+3}}|_{D^\perp},
$$
where $J$ denotes the unique, up to a sign, almost complex structure acting on the plane bundle $D^\perp$. 

It follows using \eqref{C3} that 
$$
A_{\xi_{n+2}}|_{D^\perp}\circ\Ce_T=\Ce_T^t\circ A_{\xi_{n+2}}|_{D^\perp} \,\ \text{and} \,\ A_{\xi_{n+3}}|_{D^\perp}\circ\Ce_T=\Ce_T^t\circ A_{\xi_{n+3}}|_{D^\perp}
$$
 for any $T\in\Gamma(D)$. Since $\hat{N}^f_1$ is a plane bundle, the above imply that $\Ce_T\in {\rm{span}}\{I,J\}\subseteq {\rm{End}}(D^\perp)$.
This, combined with Lemma \ref{char}, yields
$$
(\rho-1)\psi_1(T)=0 \,\ \text{and} \,\ (\rho-1)\psi_2(T)=0
$$
for any $T\in \Gamma(D)$. Thus, the splitting tensor vanishes identically on $M_3^*$. \qed
\vspace{2ex}

Hereafter, we assume that $M_3$ is not dense on $M^n$ and consider the open subset
$$
M_2=\left\{x\in M^n\smallsetminus\overline{M}_3:\dim N_1^f(x)=2\right\}.
$$
In the sequel, we assume that the open subset $M_2^*:=M_2\smallsetminus\{x\in M^n: H(x)=0\}$ is nonempty. Choose a local orthonormal frame such that  $\xi_{n+1}$ and $\xi_{n+2}$ span the plane bundle $N_1^f$ on this subset and $\xi_{n+1}$ is collinear to the mean curvature vector field. Thus, there exist smooth functions $\lambda,\mu$ such that
\begin{align*}
A_{\xi_{n+2}}e_1=\lambda e_1 + \mu e_2,  \,\ A_{\xi_{n+2}}e_2=\mu e_1 - \lambda e_2 \,\,\,\, \text{and}\,\,\,\, \lambda^2+\mu^2>0.
\end{align*}

We proceed with some auxiliary lemmas.

\begin{lemma}\po\label{n+2}
The plane bundle $N_1^f$ is parallel  in the normal connection along the distribution $D$ on the subset $M_2^*$. Moreover, the following holds: 
\begin{align}
&\mu \psi_1(T)=-\lambda \psi_2(T), \label{Nostr1}\\
&\mu \phi(T)=-(\lambda^2+\mu^2)\frac{\rho-1}{nH}\psi_2(T),\label{Nostr2}\\
&T(\mu)+2\lambda\omega(T)+(\rho+1)\lambda \psi_2(T)=0, \label{rm}\\
&T(\lambda)-2\mu\omega(T)-\mu \rho \psi_2(T)-\lambda \psi_1(T)=\frac{n\rho H}{\rho-1}\phi(T), \label{rl1}\\ 
&T(\lambda)-2\mu\omega(T)-\mu \psi_2(T)-\lambda\rho \psi_1(T)=\frac{nH}{\rho-1}\phi(T)  \label{rl2}
\end{align}
for any $T\in \Gamma(D)$, where $\phi$ is the normal connection form given by $\phi=\<\nabla^\perp\xi_{n+1},\xi_{n+2}\>.$ 
\end{lemma}

\proof
It follows from (\ref{cod}) that 
$$
\<\nabla^\perp_T\xi_{\alpha}, \xi\>=0 \;\;\text{if}\;\; \alpha=n+1,n+2
$$
for any $T\in \Gamma(D)$ and any $\xi \in \Gamma({N^f_1}^{\perp})$. Thus, the subbundle $N_1^f$ is parallel  in the normal connection along the distribution $D.$ 

Moreover, from \eqref{cod} we have 
$$
(\n_TA_{\xi_{n+2}})e_i=A_{\xi_{n+2}}\circ\Ce_Te_i+ A_{\n^{\perp}_T\xi_{n+2}}e_i, \;\;i=1,2,
$$
for any $T\in \Gamma(D)$. Bearing in mind the form of the splitting tensor given in Lemma \ref{char}, the above equations yield directly \eqref{rl1}, \eqref{rl2} and the following 
\begin{align*}
&T(\mu)+2\lambda\omega(T)+\lambda\rho \psi_2(T)-\mu \psi_1(T)=0, \\
&T(\mu)+2\lambda\omega(T)-\mu\rho \psi_1(T)+\lambda \psi_2(T)=0
\end{align*}
for any $T\in \Gamma(D)$. 
Subtracting the above equations, we obtain (\ref{Nostr1}). Similarly, \eqref{Nostr2} follows by subtracting  \eqref{rl1}, \eqref{rl2} and using  (\ref{Nostr1}).
Finally, plugging (\ref{Nostr1}) into the first of the above equations, we obtain (\ref{rm}).  \qed
\vspace{2ex}
 
Now suppose that the subset ${M_3\cup M_2}$ is not dense on $M^n$ and consider the open subset
$$
M_1=\left\{x\in M^n\smallsetminus \overline{M_3\cup M_2}: \dim N_1^f(x)=1\right\}.
$$ 

\begin{lemma}\po\label{M1}
If the subset $M_1^*:=M_1\smallsetminus\{x\in M^n: H(x)=0\}$ is nonempty, then $c=0$ and $f|_{M_1^*}$ is locally a cylinder either over a surface in $\R^{p+2}$ or over a curve in $\R^{p+1}$.
\end{lemma}

\proof
On the subset $M_1^*$ we choose a local orthonormal frame $\xi_{n+1},\dots,\xi_{n+p}$  in the normal bundle such that $\xi_{n+1}$ is collinear to the mean curvature vector field. Then we have 
$A_{\xi_\alpha}=0$ for all $\alpha\geq n+2.$ 
The Codazzi equation yields 
$$
A_{\nap_{e_i}\xi_\alpha}e_r=A_{\nap_{e_r}\xi_\alpha}e_i
$$
for all $\alpha\geq n+2,$ $i=1,2,$ and $r\geq3$. Thus, we obtain
$\nap_{e_r}\xi_{n+1}=0$ and Lemma \ref{eq} gives
\begin{equation}
2\rho(u_ru_s + v_rv_s)=c\delta_{rs}
\label{first2}
\end{equation}
for all $r\geq 3$.
Moreover, \eqref{ero} becomes
\begin{equation*}
e_r(\rho)=-\rho(\rho-1)u_r.
\end{equation*}
Differentiating (\ref{first2}) with respect to $e_r$ and using the above along with \eqref{ru} and \eqref{rv2}, we obtain
$$
\rho u_r(\rho-3)(u_r^2+v_r^2)-2c\rho u_r + 2\rho\sum_{s\geq3}^n\omega_{rs}(e_r)\big(u_su_r +v_sv_r\big)=0
$$ 
for all $r\geq3$.
In view of (\ref{first2}), the above equation simplifies to the following
$$
c(\rho+1)u_r=0.
$$

Now we prove that $c=0$. Arguing indirectly, we suppose that $c\neq0.$ 
Assume that the open set of points where $\rho\neq-1$ is nonempty. On this subset, we have $u_r=0$ for all $r\geq3.$ Thus, \eqref{first2} becomes $2\rho v_r^2=c$ for all $r\geq3$.  
Using \eqref{omega}, \eqref{rv2} yields $2\rho^2 v_r^2=c(\rho+1),$
which is a contradiction. 
Assume now that the set of points where $\rho=-1$ has nonempty interior. On this subset, \eqref{1} yields $u_r=0$ and \eqref{3} implies that $v_r=0,$ which contradicts the assumption that $c\neq0.$ 

Hence, $c=0$ and \eqref{first2} becomes
$$
\rho(u_r^2+v_r^2)=0
$$
for all $r\geq3.$ If $\rho\neq0,$ then the splitting tensor vanishes and Proposition \ref{cylinder} implies that $f$ is locally a cylinder over a surface. If the subset of points where $\rho=0$ has nonempty interior, then the Codazzi equation implies that the tangent bundle splits as an orthogonal sum of two parallel distributions one of which has rank $n-1$. Thus,  the manifold splits locally as a Riemannian product by the De Rham decomposition theorem. Since the second fundamental form is adapted to this splitting, the result follows from \cite[Th. 8.4]{da}.  \qed

\section{Submanifolds of dimension $n\geq4$}
We are now ready to give the proof of our first main result.
\vspace{2ex}

\noindent\emph{Proof of Theorem \ref{main}:}
If the open subset $M_3^*$ is nonempty, then Lemma \ref{M3} implies that the splitting tensor vanishes identically on it. Then, by Proposition \ref{cylinder} the immersion $f$ is locally a cylinder  over a surface on $M_3^*$.

Now assume that the subset $M_3$ is not dense on $M^n$ and suppose that $M_2^*$ is nonempty. Hereafter, we work on $M_2^*$.  
Due to the choice of the local orthonormal frame $\xi_{n+1}, \xi_{n+2}$ in the normal subbundle $N_1^f$, and using (\ref{Nostr1}) and  (\ref{Nostr2}), \eqref{second} of Lemma \ref{eq} takes the following form
\begin{equation}
v_rv_s\left(\lambda^2+\mu^2\right)\Big(2\rho-\left(\lambda^2+\mu^2\right)\frac{(\rho-1)^2}{n^2H^2}\Big)=c\mu^2 \delta_{rs} \label{compare}
\end{equation}
for any $r,s\geq3$. 

We claim that $v_r=0$ for any $r\geq3$. In fact, at points where 
$$
2\rho-\left(\lambda^2+\mu^2\right)\frac{(\rho-1)^2}{n^2H^2}\neq0,
$$
it follows from \eqref{compare} that
\begin{align*}
v_r^2=\frac{c\mu^2}{\left(\lambda^2+\mu^2\right)\left(2\rho-\left(\lambda^2+\mu^2\right)\frac{(\rho-1)^2}{n^2H^2}\right)}
\end{align*}
for any $r\geq3$ and
$v_rv_s=0$
for $r\neq s\geq3$.
Thus, $v_r=0$ for any $r\geq3$ at those points. 

It remains to prove that the same holds on the subset $U\subseteq M_2^*$ of points where 
\begin{equation*}
2\rho-\left(\lambda^2+\mu^2\right)\frac{(\rho-1)^2}{n^2H^2}=0.
\end{equation*}
Notice that because of \eqref{k1k2}, the subset $U$ is the set of points where 
\be\label{hypothesis}
\lambda^2+\mu^2=-2k_1k_2.
\ee
In order to prove that  $v_r=0$ for any $r\geq3$ on $U$, we suppose that the interior  of $U$ is nonempty. Suppose to  the contrary that there exists  $r_0\geq3$ such that $v_{r_0}\neq0$ on an open subset of $U$. 
Differentiating \eqref{hypothesis} with respect to $e_{r_0}$ and using \eqref{k1k2},  \eqref{1},  \eqref{2},  \eqref{rm}  and \eqref{rl1}, we obtain
$$
\lambda^2u_{r_0}-\lambda \mu v_{r_0} + (\rho+1)k_1k_2 u_{r_0}=\lambda(k_1-2k_2)\phi(e_{r_0}).
$$ 
Multiplying by $\mu$ the above and using \eqref{Nostr2}, we find that
$$
\mu u_{r_0}\left(\lambda^2+ (\rho+1)k_1k_2\right)=
\lambda v_{r_0}\Big(\mu^2- (\lambda^2+\mu^2)(k_1-2k_2)\frac{\rho-1}{nH}\Big).
$$
Taking into account \eqref{k1k2},  \eqref{Nostr1} and \eqref{hypothesis}, the above yields
$$
\lambda v_{r_0}(\rho+1)(\lambda^2+\mu^2)=0.
$$ 
Due to \eqref{hypothesis}, we conclude that $\lambda=0$ and consequently $\mu\neq0$. 
Then, it follows from \eqref{Nostr1} that $u_s=0$ for any $s\geq3$. It is easily seen from \eqref{rv1} and \eqref{rv2} for $s= r_0$  that
$$
\rho v^2_{r_0}+(\rho-1)v_{r_0}\omega(e_{r_0})-c=0
\;\;{\text {and}}\;\;
\rho v^2_{r_0}-(\rho-1)v_{r_0}\omega(e_{r_0})-c=0.
$$
Hence $\omega(e_{r_0})=0$, and consequently \eqref{omega} yields
$$
\rho v_{r_0}+\frac{\rho-1}{nH}\mu\phi(e_{r_0})=0.
$$
Using  \eqref{k1k2}, \eqref{Nostr2} and  \eqref{hypothesis}  we find that $\rho=0$, which contradicts  \eqref{hypothesis}. 
Thus we have proved the claim that $v_r=0$ for any $r\geq3$.

Now, we claim that $u_r=0$ for any $r\geq3$. It follows using \eqref{Nostr1} that $
\mu u_r=0
$
for any $r\geq3$. Obviously, the function $u_r$ vanishes at points where $\mu\neq0$. 

Assume that the set of points where $\mu=0$ has nonempty interior and argue on this subset. Since $\lambda\neq0$ on this subset, it follows from \eqref{rm} that $\2=0$ for any $r\geq3,$ and consequently \eqref{rl1} and \eqref{rl2} yield
\be
\1=\frac{\rho-1}{n H}\lambda u_r\,\,\, \text{and} \,\,\,
e_r(\lambda)=(\rho+1)\lambda u_r\label{28}
\ee
for all $r\geq3$. 
Using the first of the above equations, (\ref{second}) is written equivalently as
\begin{align}\label{c=0}
u_r^2\Big(2\rho-\frac{\lambda^2(\rho-1)^2}{n^2H^2}\Big)=c
\end{align} 
for all $r\geq3$.

Since we already proved that $v_r=0$ for all $r\geq3$, Lemma \ref{char} implies that the image of the splitting tensor $C\colon D\to {\rm{End}}(D^\perp)$ satisfies $\dim {\rm{Im}}\, C\leq1.$
Thus, $\dim\ker\, C\geq n-3.$

Now suppose that $\dim\ker\, C=n-3$. Then, there  exists a unique $r_0\geq3$ such that $u_{r_0}\neq0$ and $u_s=0$ for any $s\neq r_0$. Thus, \eqref{c=0} implies that $c=0$ and 
$$
2\rho=\frac{\lambda^2(\rho-1)^2}{n^2H^2}.
$$
On account of \eqref{k1k2}, the above equation becomes $\lambda^2=-2k_1k_2>0.$
Differentiating this equation with respect to $e_{r_0}$ and using \eqref{k1k2},  \eqref{1}, \eqref{2} and the second of \eqref{28}, we obtain $2\lambda^2 + k_1k_2=0,$
which contradicts the previous equation. Thus, the splitting tensor vanishes identically on the subset $M_2^*$ and consequently, by Proposition \ref{cylinder}, the immersion $f$ is locally a cylinder over a surface.

If the open subset $M_1^*$ is nonempty, then Lemma \ref{M1} implies that $f$ is locally a cylinder over a surface or over a curve. \qed

\section{Elliptic submanifolds}

In this section, we recall from \cite{DF2} the notion of elliptic submanifolds 
of a space form as well as several of their basic properties.

Let $f\colon M^n\to\Q_c^m$ be an isometric immersion.
The $\ell^{th}$\emph{-normal space} $N^f_\ell(x)$  of $f$
at $x\in M^n$ for $\ell\ge 1$ is defined as
$$
N^f_\ell(x)=\spa\left\{\alpha^f_{\ell+1}\left(X_1,\ldots,X_{\ell+1}\right):
X_1,\ldots,X_{\ell+1}\in T_xM\right\}.
$$
Here $\alpha^f_2=\alpha^f$ and for $s\geq 3$ the so-called 
$s^{th}$\emph{-fundamental form} is the symmetric tensor 
$\alpha^f_s\colon TM\times\cdots\times TM\to N_fM$  defined inductively by
$$
\alpha^f_s(X_1,\ldots,X_s)=\pi^{s-1}\left(\nabla^\perp_{X_s}\cdots
\nabla^\perp_{X_3}\alpha^f(X_2,X_1)\right),
$$
where $\pi^k$ stands for the projection onto $(N_1^f\oplus\cdots\oplus N_{k-1}^f)^{\perp}$.

An isometric immersion $f\colon M^n\to\Q_c^m$ is called \emph{elliptic}  
if $M^n$ carries a totally geodesic distribution $D$ of rank $n-2$ satisfying $D(x)\subseteq \D_f(x)$ for any $x\in M^n$ and  there exists an (necessary unique up to a sign) almost complex 
structure $J\colon D^\perp\to D^\perp$ such that  the second fundamental form satisfies
$$
\alpha^f(X,X)+\alpha^f(JX,JX)=0
$$
for all $X\in D^\perp$. Notice that $J$ is orthogonal if and only $f$ is 
minimal. 

Assume that $f\colon M^n\to\Q_c^m$ is  substantial  and elliptic.
Assume also that $f$ is \emph{nicely curved} which means that for any $\ell\geq 1$
all subspaces $N^f_\ell(x)$ have constant dimension and thus form subbundles
of the normal bundle. Notice that any $f$ is nicely curved along connected 
components of an open dense subset of $M^n$.
Then, along that subset the normal bundle splits orthogonally and smoothly as
\be\label{splits}
N_fM=N^f_1\oplus \cdots \oplus N^f_{\tau_f},
\ee
where all $N^f_\ell$'s have rank two, except possibly the last one 
that has rank one in case the codimension is odd.
Thus, the induced bundle $f^*T\Q_c^m$ splits as 
$$
f^*T\Q_c^m=f_*D\oplus N^f_0\oplus N^f_1\oplus \cdots \oplus N^f_{\tau_f},
$$
where $N^f_0=f_*D^\perp$. Setting
$$
\tau^o_f = \left\{\begin{array}{l}
\tau_f\;\;\;\;\;\;\;\;\;\;\;\mbox{if}\;\;m-n\;\;\; \mbox{is even}\\
\tau_f-1\;\;\;\;\; \mbox{if}\;\;m-n\;\;\; \mbox{is odd}
\end{array} \right.
$$ 
it turns out that the almost complex structure $J$ on $D^\perp$ induces an almost
complex  structure $J_\ell$ on each $N_\ell^f$, $0\leq \ell\leq\tau^o_f$, defined by
$$
J_\ell\alpha_{\ell+1}^f\left(X_1,\ldots,X_\ell,X_{\ell+1}\right)
=\alpha_{\ell+1}^f\left(X_1,\ldots,X_\ell,J X_{\ell+1}\right),
$$
where $\alpha_1^f=f_*$.

The \emph{$\ell^{th}$-order curvature ellipse} 
$\mathcal{E}_\ell^f(x)\subset N^f_\ell(x)$ of $f$ at $x\in M^n$ for $0\leq\ell\leq\tau^o_f$ is 
$$
\mathcal{E}_\ell^f(x)=\big\{\alpha^f_{\ell+1}(Z_{\theta},\dots,Z_{\theta}): 
Z_{\theta}=\cos\theta Z+\sin\theta J Z\;\;\mbox{and}\;\;\theta\in [0,\pi)
\big\},
$$
where $Z\in D^\perp(x)$ has  unit length and satisfies $\<Z,JZ\>=0$. From 
ellipticity such a $Z$ always exists and $\mathcal{E}_\ell^f(x)$ is indeed  
an ellipse. 

We say that the curvature ellipse $\mathcal{E}_\ell^f$ of an elliptic submanifold
$f$ is a \emph{circle} for some $0\leq\ell\leq\tau^o_f$ if all ellipses 
$\mathcal{E}_\ell^f(x)$ are circles. That the curvature ellipse 
$\mathcal{E}_\ell^f$ is a circle is equivalent to the almost complex  
structure $J_\ell$ being orthogonal.
Notice that $\mathcal{E}_0^f$ is a circle if and only if $f$ is minimal.

Let  $f\colon M^n\to\Q_c^{m-c}, c\in \{0,1\},$ be a substantial nicely curved 
elliptic submanifold. Assume that $M^n$ is the saturation  of a fixed cross 
section $L^2\subset M^n$ to the foliation of the distribution $D.$
The subbundles in the orthogonal splitting (\ref{splits}) are parallel in 
the normal connection (and thus in $\Q_c^{m-c}$) along $D$. 
Hence each $N^f_\ell$ can be seen as a vector bundle along the 
surface $L^2$.

A \emph{polar surface} to $f$ is an immersion $h$ of $L^2$
defined as follows:
\begin{itemize}
\item [(a)] If $m-n-c$ is odd, then  the polar surface $h\colon L^2\to\Sf^{m-1}$ 
is the spherical image of the unit normal field spanning $N^f_{\tau_f}$. 
\item [(b)] If $m-n-c$ is even, then  the polar surface $h\colon L^2\to\R^m$ is 
any surface such that $h_*T_xL=N^f_{\tau_f}(x)$ up to parallel 
identification in $\R^m$.
\end{itemize}

Polar surfaces always exist since in case ${\rm(b)}$  any elliptic submanifold 
admits locally many polar surfaces. 

The almost complex structure $J$ on  $D^\perp$ induces an almost complex 
structure $\tilde J$ on $TL$ defined by $P\circ\tilde J=J\circ P$,
where $P\colon TL \to D^\perp$ is the orthogonal projection.
It turns out that a polar surface to an elliptic submanifold is necessarily elliptic.  
Moreover, if the elliptic submanifold has a circular curvature ellipse  then its 
polar surface has the same property at the ``corresponding" normal bundle.
As a matter of fact, up to parallel identification it holds that
\be\label{eqp}
N_s^h=N_{\tau^o_f-s}^f\;\;\mbox{and}\;\;
J^h_s=\big(J^f_{\tau^o_f-s}\big)^t,\;\; 0\leq s\leq\tau^o_f.
\ee
In particular, the polar surface is nicely curved.

A \emph{bipolar surface} to $f$  is any polar surface to a polar surface to $f$.
In particular, if we are in case $f\colon M^3\to\Sf^{m-1}$, then a bipolar
surface to $f$ is a nicely curved elliptic surface $g\colon
L^2\to\R^m$.

\section{Three-dimensional submanifolds}

In this section, we study the case of three-dimensional submanifolds and we provide the proof of Theorem \ref{main1}.
To this purpose, we need the following results.

\begin{proposition}\po\label{flat}
Let $f\colon M^3 \to \mathbb{Q}^{3+p}_c$ be an isometric immersion such that $M^3$ carries a totally geodesic distribution $D$ of rank one satisfying $D(x)\subseteq \D_f(x)$ for any $x\in M^3$. If the mean curvature of $f$ is constant along each integral curve of $D$ and the normal bundle of $f$ is flat, then $f$ is minimal or $c=0 $ and $f$ is locally a cylinder.
\end{proposition}
\proof
Assume that $f$ is nonminimal.  If the open subset $M_3^*$ is nonempty, then Lemma \ref{M3} and Proposition \ref{cylinder} imply that the immersion $f$ is locally a cylinder over a surface.

Now suppose that the open subset $M_2^*$ is nonempty and argue on it.
Having flat normal bundle implies that $\mu=0$ and  according to \eqref{Nostr1} we obtain $v_3=0$. Consequently, \eqref{rv2} is written as
\be\label{a3}
e_3(u_3)=u_3^2+c.
\ee 
Comparing \eqref{rl1} and \eqref{rl2}, we obtain
$$
\phi(e_3)=\frac{\rho-1}{nH}\lambda u_3.
$$
Thus, 
\be\label{e3l}
e_3(\lambda)=(\rho+1)\lambda u_3
\ee
and consequently \eqref{ero} becomes
\be\label{ero'}
e_3(\rho)=u_3(\rho-1)(\tau-\rho),
\ee
where $\tau $ is the function given by
$$
\tau=\frac{\lambda^2(\rho-1)^2}{n^2H^2}.
$$
Moreover, \eqref{second} is written as $u_3^2(2\rho-\tau)=c.$
Differentiating with respect to $e_3$ and using \eqref{a3}-\eqref{ero'}, we derive that
$$
u_3^2(\rho+1)(\rho-\tau)=0.
$$

Now we claim that $u_3=0.$ Arguing indirectly, we suppose that $u_3\neq0$ on an open subset. Observe that 
$\rho\neq-1$ due to our assumption and \eqref{1}. Hence $\rho=\tau,$  or  equivalently $\rho n^2H^2=\lambda^2(\rho-1)^2$ and $e_3(\rho)=0$ by \eqref{ero'}. Thus $e_3(\lambda)=0$, which contradicts \eqref{e3l}  since $\lambda\neq0.$ This proves the claim that $u_3=0$ and consequently the splitting tensor vanishes. That the immersion $f$ is locally a cylinder on $M_2^*$ follows from Proposition \ref{cylinder}.

If the open subset $M_1^*$ is nonempty, then Lemma \ref{M1} implies that the immersion $f$ is locally a cylinder over a  surface or over a curve. \qed

\begin{proposition}\po \label{sphelliptic}
Let $f\colon M^3 \to \mathbb{Q}^{3+p}_c$ be a nonminimal isometric immersion such that $M^3$ carries a totally geodesic distribution $D$ of rank one satisfying $D(x)\subseteq \D_f(x)$ for any $x\in M^3$. If the mean curvature of $f$ is constant along each integral curve of $D$ and $f$ is not locally a cylinder, then the splitting tensor of $f$ is an almost complex structure on $D^\perp.$ Moreover, $f$ is a spherical elliptic submanifold with respect to this almost complex structure and its first curvature ellipse is a circle.
\end{proposition}

\proof
Since by assumption the immersion $f$ is not a cylinder on any open subset, it follows from Proposition \ref{cylinder}, Lemmas \ref{M3} and  \ref{M1} that the open subsets $M_3^*$ and $M_1^*$ are both empty.

Proposition \ref{flat}, implies that the immersion $f$ has nonflat normal bundle on $M_2^*.$
Thus, we have $\mu\neq0$ and $\rho\neq-1.$ Using \eqref{Nostr1} and \eqref{Nostr2}, it is easily seen that  \eqref{second},  \eqref{ero},  \eqref{ru}, \eqref{omega}, \eqref{rm} and \eqref{rl1} are written as
\begin{align}
&\omega(e_3)=-\frac{\rho-\tau}{\rho+1}v_3,\nonumber\\
&e_3(\rho)=\frac{\lambda}{\mu}(\rho-1)(\rho-\tau)v_3,\label{e3ro}\\
&e_3(\mu)=-\frac{\lambda}{\rho+1}\left(2\tau+\rho^2+1\right)v_3,\nonumber\\
&e_3(\lambda)=\Big(\frac{2\mu}{\rho+1}\tau-\frac{2\mu\rho}{\rho+1}-\frac{\lambda^2}{\mu}(\rho+1)\Big)v_3,\label{3la}\\
&e_3(v_3)=\frac{\lambda}{\mu(\rho+1)}\left((\rho-1)\tau -(2\rho^2+\rho+1)\right)v_3^2,\nonumber \\
&(\lambda^2+\mu^2)(2\rho-\tau)v_3^2=c\mu^2, \label{comp1}
\end{align}
where $\tau$ is the function given by 
$$
\tau=(\lambda^2+\mu^2)\frac{(\rho-1)^2}{n^2H^2}.
$$
By differentiating \eqref{comp1} and using all the above equations, we obtain
$$
\lambda(\lambda^2+\mu^2)\Big(\rho(5\rho^2+6\rho+5)-(4\rho^2+2\rho+4)\tau - 2\tau^2\Big)v_3^3=c\lambda\mu^2v_3.
$$

We claim that $\lambda v_3=0.$ Arguing indirectly, we assume that the open subset where $\lambda v_3\neq0$ is nonempty. Thus, comparing the above equation with \eqref{comp1}, we derive that $\tau=\rho.$ This along with \eqref{e3ro} imply that $e_3(\tau)=e_3(\rho)=0.$ By the definition of $\tau,$ it follows that $e_3(\lambda^2+\mu^2)=0.$ Using the above equations it is easy to see that 
$$
e_3(\lambda^2+\mu^2)=-2\frac{\lambda}{\mu}(\lambda^2+\mu^2)(\rho+1)v_3,
$$
which is a contradiction and this proves our claim.

Now we claim that $v_3$ cannot vanish on any open subset. Arguing indirectly, we suppose that $v_3=0$ on an open subset. Then \eqref{Nostr1}  implies that $u_3=0. $ By Lemma \ref{char}, the splitting tensor vanishes and consequently the immersion $f$ would be a cylinder by Proposition \ref{cylinder}. This contradicts our assumption.
 
Since we already proved that $\lambda v_3=0$, we obtain $\lambda=0$ and \eqref{Nostr1} implies that $u_3=0.$ It follows from \eqref{3la} that
\be\label{romi}
\mu^2=\frac{\rho n^2H^2}{(\rho-1)^2}.
\ee
In particular, we have $\rho>0.$ This, along with \eqref{comp1} yield 
\be\label{roci}
\rho v_3^2=c.
\ee
Hence, $c=1.$ Now observe that the splitting tensor satisfies $\Ce_3^2=-I,$
where $I$ is the identity endomorphism of $D^{\perp},$ that is, $\Ce_3$ is an almost complex structure $J\colon D^\perp\to D^\perp.$ Using \eqref{roci} and the fact that the shape operator $A_{\xi_5}$ satisfies
$A_{\xi_5}e_i=\mu e_j$ for $i\neq j=1,2,$ we easily verify that the second fundamental form of $f$ satisfies
$\alpha^f(Je_1,e_2)=\alpha^f(e_1,Je_2).$
This is equivalent to the ellipticity of the immersion $f.$ 

In order to prove that the first curvature ellipse of $f$ is a circle, it is equivalent to prove that the vector fields $\alpha^f(e_1,e_1)$ and $\alpha^f(e_1,Je_1)$ are of the same length and perpendicular. Obviously, they are perpendicular since 
$$
\alpha^f(e_1,e_1)=k_1\xi_4  \;\; \text{and}  \;\;
\alpha^f(e_1,Je_1)=\mu v_3\xi_5.
$$
Using \eqref{k1k2} and  \eqref{romi}, we obtain 
$$
\frac{\|\alpha^f(e_1,Je_1)\|^2}{\|\alpha^f(e_1,e_1)\|^2}=\rho v_3^2.
$$
Bearing in mind \eqref{roci}, we conclude that the first curvature ellipse is a circle. \qed
\vspace{2ex}

The following result parametrizes all three-dimensional submanifolds in spheres that carry a totally geodesic distribution of rank one, contained in the relative nullity distribution, such that the mean curvature is constant along each integral curve. This parametrization, given in terms of their polar surfaces, was introduced in \cite{DF2} as the \textit{polar parametrization}. 

\begin{theorem}\po\label{main0}
Let $h\colon L^2\to \mathbb{Q}_c^{N+1}, c\in\{0,1\},N\geq 5,$ be a nicely curved elliptic surface of substantial even codimension, such that the curvature ellipses $\mathcal{E}_{\tau_h-2}^h, \mathcal{E}_{\tau_h}^h$  are circles and $\mathcal{E}_{\tau_h-1}^h$ is nowhere a circle. Then, the map $\Psi_h\colon M^3\to \mathbb{S}^{N+c}$ defined on the circle bundle $M^3=UN^h_{\tau_h}=\{(x,w)\in N^h_{\tau_h}: \|w\|=1\}$ by $\Psi_h(x,w)=w$ is a nonminimal elliptic isometric immersion with polar surface $h$. Moreover, $M^3$ carries a totally geodesic distribution $D$  of rank one satisfying $D(p)\subseteq\D_{\Psi_h}(p)$ for any $p\in M^3$ such that the mean curvature of $\Psi_h$ is constant along each integral curve of $D.$ 

Conversely, let $f\colon M^3 \to \mathbb{S}^{3+p}, p\geq2,$ be a substantial nonminimal isometric immersion such that $M^3$ carries a totally geodesic distribution $D$ of rank one satisfying $D(x)\subseteq \D_f(x)$ for any $x\in M^3$. If the mean curvature of $f$ is constant along each integral curve of $D$, then f is elliptic and there exists an open dense subset of $M^3$ such that 
for each point there exist a neighborhood $U,$  and a local isometry $F\colon U\to UN^h_{\tau_h}$ such that $f=\Psi_h\circ F,$ where  $h$  is a polar surface to $f$ with curvature ellipses as above.
\end{theorem}

\proof
Let $h\colon L^2\to \mathbb{Q}_c^{N+1}, c\in\{0,1\},$ be a substantial elliptic surface, where $N=2m+3, m\geq1$.
Choose a local orthonormal frame $e_1, e_2$ in the tangent bundle of $L^2$ such that the almost complex structure $J$ of the elliptic surface is given by
$$
Je_1=be_2 \,\ \, \text{and} \,\ \, Je_2=-\frac{1}{b}e_1,
$$
where $b$ is a positive smooth function.

We argue for the case where $m\geq2.$ The case where $m=1$ can be handled in a similar manner.
We know from \eqref{splits} that the normal bundle splits orthogonally as 
$$
N_hL=N^h_1\oplus\cdots\oplus N^h_{m-1}\oplus N^h_m\oplus N^h_{m+1}.
$$
Let $\zeta_3,\dots,\zeta_{2m+4}$ be an orthonormal frame  in the normal bundle, defined on an open subset $V\subseteq L^2,$ such that $\zeta_{2s+1}, \zeta_{2s+2}$ span the plane subbundle $N^h_s$ for any $1\leq s\leq m+1.$ The corresponding normal connection forms $\omega_{\alpha\beta}$ are given by $\omega_{\alpha\beta}=\<\nap \zeta_\alpha,\zeta_\beta\>, \alpha,\beta=3,\dots,2m+4.$

Due to our hypothesis, we may choose the frame such that 
$$
\alpha^h_m(e_1,\dots,e_1)=\kappa_{m-1}\zeta_{2m-1}, \,\ \,\  \alpha^h_m(e_1,\dots,e_1,e_2)=\frac{\kappa_{m-1}}{b}\zeta_{2m}
$$
and
$$
\alpha^h_{m+2}(e_1,\dots,e_1)=\kappa_{m+1}\zeta_{2m+3}, \,\ \,\ \alpha^h_{m+2}(e_1,\dots,e_1,e_2)=\frac{\kappa_{m+1}}{b}\zeta_{2m+4},
$$
where $\kappa_{m-1}, \kappa_{m+1}$ denote the radii of the circular curvature ellipses $\mathcal{E}^h_{m-1}, \mathcal{E}^h_{m+1},$ respectively. Since the curvature ellipse $\mathcal{E}^h_{m}$ is nowhere a circle, we may choose $\zeta_{2m+1}, \zeta_{2m+2}$ to be collinear to the major and minor axes of this ellipse, respectively. Thus, we may write 
$$
\alpha^h_{m+1}(e_1,\dots,e_1)=v_{11}\zeta_{2m+1}+v_{12}\zeta_{2m+2} \,\ \,\ \text{and} \,\ \,\ \alpha^h_{m+1}(e_1,\dots,e_1,e_2)=v_{21}\zeta_{2m+1}+v_{22}\zeta_{2m+2},
$$
where $v_{ij}$ are smooth functions such that 
\be\label{axes}
b^2v_{21}v_{22}+v_{11}v_{12}=0, \,\ \,\ \kappa_m=\left(v_{11}^2+b^2v_{21}^2\right)^{1/2},  \,\ \,\ \mu_{m}=\left(v_{12}^2+b^2v_{22}^2\right)^{1/2}
\ee
and $\kappa_m,\mu_m$ denote the lengths of the semi-axes of the curvature ellipse $\mathcal{E}^h_{m}$.

Bearing in mind the definition of the higher fundamental forms, their symmetry and the ellipticity of the surface $h,$ we have
$$
\alpha_{s+1}^h\left(e_1,\dots,e_1,e_2\right)=\left(\nap_{e_2}\alpha_s^h\left(e_1,\dots,e_1\right)\right)^{N^h_s}=\left(\nap_{e_1}\alpha_s^h\left(e_1,\dots,,e_1,e_2\right)\right)^{N^h_s},
$$
$$
\alpha_{s+1}^h\left(e_1,\dots,e_1\right)=-b^2\left(\nap_{e_2}\alpha_s^h\left(e_1,\dots,e_1, e_2\right)\right)^{N^h_s}=\left(\nap_{e_1}\alpha_s^h\left(e_1,\dots,,e_1\right)\right)^{N^h_s}
$$
for $s=m,m+1,$ where $( \cdot )^{N^h_s}$ denotes taking the projection onto the normal subbundle $N^h_s$. 
From these we obtain
\begin{align}
&\omega_{2m-1,2m+1}(e_1)=\frac{v_{11}}{\kappa_{m-1}},\label{-11} \,\ \,\ \omega_{2m-1,2m+2}(e_1)=\frac{v_{12}}{\kappa_{m-1}},\\
&\omega_{2m-1,2m+1}(e_2)=\frac{v_{21}}{\kappa_{m-1}},\label{-12} \,\ \,\ \omega_{2m-1,2m+2}(e_2)=\frac{v_{22}}{\kappa_{m-1}},\\
&\omega_{2m,2m+1}(e_1)=\frac{b v_{21}}{\kappa_{m-1}},\label{01} \,\ \,\ \omega_{2m,2m+2}(e_1)=\frac{b v_{22}}{\kappa_{m-1}},\\
&\omega_{2m,2m+1}(e_2)=-\frac{v_{11}}{b \kappa_{m-1}},\label{02} \,\ \,\
\omega_{2m,2m+2}(e_2)=-\frac{v_{12}}{b \kappa_{m-1}},  \\
&\omega_{2m+1,2m+3}(e_1)=\frac{b \kappa_{m+1}}{\kappa_m \mu_m}v_{22},\label{131} \,\ \,\ \omega_{2m+1,2m+3}(e_2)=\frac{\kappa_{m+1}}{b \kappa_m \mu_m}v_{12}, \\
&\omega_{2m+1,2m+4}(e_1)=-\frac{\kappa_{m+1}}{\kappa_m \mu_m}v_{12},\label{141} \,\ \,\ \omega_{2m+1,2m+4}(e_2)=\frac{\kappa_{m+1}}{\kappa_m \mu_m}v_{22},\\
&\omega_{2m+2,2m+3}(e_1)=-\frac{b \kappa_{m+1}}{\kappa_m \mu_m}v_{21}, \,\ \,\ \omega_{2m+2,2m+3}(e_2)=-\frac{\kappa_{m+1}}{b \kappa_m \mu_m}v_{11},\label{.}\\
&\omega_{2m+2,2m+4}(e_1)=\frac{\kappa_{m+1}}{\kappa_m \mu_m}v_{11}, \,\ \,\
\omega_{2m+2,2m+4}(e_2)=-\frac{\kappa_{m+1}}{\kappa_m \mu_m}v_{21}.  \label{..} 
\end{align}

Let $\Pi\colon M^3\to L^2$ the natural projection of the circle bundle 
$$
M^3=UN^h_{\tau_h}=\left\{(x,\delta)\in N^h_{m+1}\colon \|\delta\|=1, x\in L^2\right\}.
$$
We parametrize $\Pi^{-1}(V)$ by $V\times\R$ via the map 
$$
(x,\theta)\mapsto\left(x,\cos\theta \zeta_{2m+3}(x)+\sin\theta \zeta_{2m+4}(x)\right)
$$
and consequently, we have
$$
\Psi_h(x,\theta)=\cos \theta \zeta_{2m+3}+\sin \theta \zeta_{2m+4}.
$$ 

Notice that $\nap N^h_{m+1}\subseteq N_m^h\oplus N_{m+1}^h.$ It is easily seen that
\begin{align*}
\Psi_{h_*}E_i&=\big(\cos\theta\omega_{2m+3,2m+1}(e_i)+\sin\theta\omega_{2m+4,2m+1}(e_i)\big)\zeta_{2m+1}\\
&+\big(\cos\theta\omega_{2m+3,2m+2}(e_i)+\sin\theta\omega_{2m+4,2m+2}(e_i)\big)\zeta_{2m+2},
\end{align*}
where the vector fields $E_i\in TM,i=1,2,$ are given by 
$$
E_i=e_i-\omega_{2m+3,2m+4}(e_i)\frac{\partial}{\partial\theta}.
$$
Using \eqref{131}-\eqref{..}, we obtain
\be\label{dPsi}
\Psi_{h_*}E_1=\frac{\kappa_{m+1}}{\kappa_m\mu_m}\Big(\big(- b v_{22}\cos\theta+v_{12}\sin\theta \big)\zeta_{2m+1}
+\big(b v_{21}\cos\theta - v_{11}\sin\theta \big)\zeta_{2m+2}\Big)
\ee
and
\be\label{dPsi'}
\Psi_{h_*}E_2=\frac{\kappa_{m+1}}{\kappa_m\mu_m}\Big(-\big(\frac{v_{12}}{b}\cos\theta+ v_{22}\sin\theta\big)\zeta_{2m+1}+\big(\frac{v_{11}}{b}\cos\theta +  v_{21}\sin\theta\big)\zeta_{2m+2}\Big).
\ee
Additionally, we have
\be\label{theta}
\Psi_{h_*}(\partial/\partial\theta)=-\sin\theta \zeta_{2m+3}+\cos\theta \zeta_{2m+4}.
\ee

It  follows that the normal bundle of the isometric immersion $\Psi_h$ is given by
$$
N_{\Psi_h}M=c\ {\rm{span}} \{h\}\oplus N^h_1\oplus\cdots\oplus N^h_{m-2}\oplus N^h_{m-1}.
$$
It is easy to see that the first normal bundle of $\Psi_h$ is
$
N_1^{\Psi_h}=N^h_{m-1}.
$
Moreover, it follows easily that the distribution $D={\rm{span}}\{\partial/\partial\theta\}$ is contained in the nullity distribution $\D_{\Psi_h}$ of $\Psi_h.$ In particular, from \eqref{theta} and the Gauss formula we derive that 
$
\n_{\partial/\partial\theta}\partial/\partial\theta=0.
$
This implies that  the distribution $D$  is totally geodesic.

It remains to show that the mean curvature of the immersion $\Psi_h$ is constant along each integral curve of $D$.
The shape operator $A_{\zeta_{2m-j}}$ of $\Psi_h$  with respect to the normal direction $\zeta_{2m-j}, j=0,1,$ is given  by the Weingarten formula as
\be\label{-dPsi}
-\Psi_{h_*}\left(A_{\zeta_{2m-j}}E_i\right)=\nap_{e_i}\zeta_{2m-j}-\Big(\tilde{\n}_{e_i}\zeta_{2m-j}\Big)^{N^h_{m-2}\oplus N^h_{m-1}} =\left(\nap_{e_i}\zeta_{2m-j}\right)^{N^h_m},\;\;i=1,2,
\ee
since $\zeta_{2m-1},\zeta_{2m}\in N_{m-1}^h$. Here, $\tilde{\n}$ stands for the induced connection of the induced bundle $h^*T\mathbb{Q}_c^{N+1}.$
It follows from  \eqref{-dPsi}  using \eqref{-11}-\eqref{02} that 
\begin{align}
&\Psi_{h_*}\left(A_{\zeta_{2m-1}}E_1\right)=-\frac{1}{\kappa_{m-1}}\left(v_{11}\zeta_{2m+1}+v_{12}\zeta_{2m+2}\right),\label{po}\\
&\Psi_{h_*}\left(A_{\zeta_{2m-1}}E_2\right)=-\frac{1}{\kappa_{m-1}}\left(v_{21}\zeta_{2m+1}+v_{22}\zeta_{2m+2}\right),\label{popo}\\
&\Psi_{h_*}\left(A_{\zeta_{2m}}E_1\right)=-\frac{b}{\kappa_{m-1}}\left(v_{21}\zeta_{2m+1}+v_{22}\zeta_{2m+2}\right),\label{popopo}\\
&\Psi_{h_*}\left(A_{\zeta_{2m}}E_2\right)=\frac{1}{b \kappa_{m-1}}\left(v_{11}\zeta_{2m+1}+v_{12}\zeta_{2m+2}\right).
\label{popopopo}\\ \nonumber
\end{align}

We may set
\be\label{2m-1}
A_{\zeta_{2m-1}}E_i=\lambda_{i1}E_1+\lambda_{i2}E_2 \,\ \,\ \text{and} \,\ \,\ A_{\zeta_{2m}}E_i=\gamma_{i1}E_1+\gamma_{i2}E_2, \,\ i=1,2,
\ee
where $\lambda_{ij}$ and $\gamma_{ij}$ are smooth functions on the manifold $M^3$. 
From \eqref{dPsi}, \eqref{po}, \eqref{popo} and the first one of \eqref{2m-1}, we obtain
$$
\lambda_{11}=\frac{1}{\kappa_{m-1}\kappa_{m+1}}\left(\left(v_{11}^2+v_{12}^2\right)\cos\theta+ b\left(v_{11}v_{21}+v_{12}v_{22}\right)\sin\theta\right)
$$
and
$$
\lambda_{22}=\frac{1}{\kappa_{m-1}\kappa_{m+1}}\left(- b^2\left(v_{21}^2+v_{22}^2\right)\cos\theta+ b\left(v_{11}v_{21}+v_{12}v_{22}\right)\sin\theta\right).
$$
Hence
$$
{\rm trace} A_{\zeta_{2m-1}}=\frac{1}{\kappa_{m-1}\kappa_{m+1}}\left(\left(v^2_{11}+v^2_{12}- b^2v_{21}^2-b^2v_{22}^2\right)\cos\theta+ 2b\left(v_{11}v_{21}+v_{12}v_{22}\right)\sin\theta\right).
$$

Similarly, from \eqref{dPsi'}, \eqref{popopo}, \eqref{popopopo} and the second of \eqref{2m-1}, 
we find that
$$
\gamma_{11}=\frac{1}{\kappa_{m-1}\kappa_{m+1}}\left(b\left(v_{11}v_{21}+v_{12}v_{22}\right)\cos\theta+ b^2\left(v_{21}^2+v_{22}^2\right)\sin\theta\right)
$$
and
$$
\gamma_{22}=\frac{1}{\kappa_{m-1}\kappa_{m+1}}\left(b\left(v_{11}v_{21}+v_{12}v_{22}\right)\cos\theta-\left(v_{11}^2+v_{12}^2\right)\sin\theta\right).
$$
Then, it follows that
$$
{\rm trace} A_{\zeta_{2m}}=\frac{1}{\kappa_{m-1}\kappa_{m+1}}\left(2b\left(v_{11}v_{21}+v_{12}v_{22}\right)\cos\theta-\left(v^2_{11}+v^2_{12}- b^2v_{21}^2-b^2v_{22}^2\right) \sin\theta\right).
$$
Thus, the mean curvature of the isometric immersion $\Psi_h$ is given by
$$
\|\h_{\Psi_h}\|^2=\frac{1}{(3\kappa_{m-1}\kappa_{m+1})^2}\Big(\big(v_{11}^2+v_{12}^2+b^2v_{21}^2+b^2v_{22}^2\big)^2-4\big(v_{11}^2+b^2v_{21}^2\big)^2\big(v_{12}^2+b^2v_{22}^2\big)^2\Big).
$$
Using \eqref{axes}, the above equation becomes
$$
\|\h_{\Psi_h}\|=\frac{|\kappa^2_m-\mu^2_m|}{3\kappa_{m-1}\kappa_{m+1}}.
$$
It is clear that the mean curvature of the isometric immersion $\Psi_h$ is constant along each integral curve of the distribution $D.$ This completes the proof of the direct statement of the theorem for $m\geq2.$ 
The case $m=1$ can be treated in a similar manner. In this case, the mean curvature of $\Psi_h$ is given by
$$
\|\h_{\Psi_h}\|=\frac{|\kappa^2_1-\mu^2_1|}{3\kappa_{2}^2}.
$$

Conversely,  let $f\colon M^3\to\mathbb{S}^{3+p}$ be a nonminimal isometric immersion. Suppose that $M^3$ carries a totally geodesic distribution $D$ of rank one satisfying $D(x)\subseteq \D_f(x)$ for any $x\in M^3$ such that the mean curvature is constant along each integral curve of $D.$ 
From Proposition \ref{sphelliptic}, we know that $f$ is an elliptic submanifold and its first curvature ellipse is a circle. Hereafter, we work on a connected component of an open dense subset where $f$ is nicely curved.

Consider a polar surface $h\colon L^2\to \mathbb{Q}^{p-c+4}_c$  to the immersion $f,$  where $c=0$ if $p$ is even and $c=1$ if $p$ is odd. Notice that $\tau^0_f=\tau_h-1.$
Using \eqref{eqp}, we conclude that the  curvature ellipse $\mathcal{E}^h_{\tau_h-2}$ of the surface $h$ is a circle
and the curvature ellipse $\mathcal{E}^h_{\tau_h-1}$ is nowhere a circle.

We claim that the last curvature ellipse $\mathcal{E}^h_{\tau_h}$ is a circle. Observe that 
$N^h_{\tau_h}={\rm{span}}\{\xi,\eta \},$
where the sections $\xi,\eta$ of the normal bundle $N_hL$ are given by 
$\xi=f\circ \pi\;\;\text{and}\;\; \eta=f_*e_3\circ\pi.$ 
Here $\pi$  denotes the natural projection $\pi\colon M^3\to L^2$  onto the fixed cross section $L^2\subset M^3$ to the foliation generated by the distribution $D.$

Let $X_1,\dots,X_{\tau_h}\in TL$ be arbitrary vector fields. By \eqref{eqp} we have $N^h_{\tau_h-1}=N^f_0=f_*D^\perp.$ Thus, there exists $X\in \Gamma(D^\perp)$ such that 
$$
\alpha^h_{\tau_h}\left(X_1,\dots,X_{\tau_h}\right)=f_*X.
$$
For every vector field $Y\in TL$ there exists a vector field $Z\in \Gamma(D^\perp)$ such that $Y=\pi_*Z.$ Then 
we have
\begin{align*}
\alpha^h_{\tau_h+1}(X_1,\dots,X_{\tau_h},Y)&=\left(\nap_Y\alpha^h_{\tau_h}(X_1,\dots,X_{\tau_h})\right)^{N^h_{\tau_h}}\nonumber\\
&= -\<f_*X,f_*Z\>\xi-\<f_*X, \tilde{\n}_Zf_*e_3\>\eta.
\end{align*}
Using the Gauss formula and the definition of the splitting tensor, the above equation becomes
$$
\alpha^h_{\tau_h+1}(X_1,\dots,X_{\tau_h},Y)=-\<X,Z\>\xi+\<X,\Ce_3Z\>\eta.
$$ 
From Proposition \ref{sphelliptic}, we know that the splitting tensor in the direction of $e_3$ is the almost complex structure $J_0^f\colon D^\perp\to D^\perp$ of $f.$ Hence, we obtain 
$$
\alpha^h_{\tau_h+1}(X_1,\dots,X_{\tau_h},Y)=-\<X,Z\>\xi+\<X,J^f_0Z\>\eta.
$$ 
On account of $\pi_*\circ J^f_0=J^h_0\circ \pi_*,$
we have $J^h_0Y=\pi_*J^f_0Z.$
Thus, it follows that
$$
\alpha^h_{\tau_h+1}(X_1,\dots,X_{\tau_h},J^h_0Y)=-\<X,J_0^hZ\>\xi-\<X,Z\>\eta.
$$ 
Since $\xi,\eta $ is an orthonormal frame of the subbundle $N^h_{\tau_h},$ it is now obvious that the normal vector fields $\alpha^h_{\tau_h+1}(X_1,\dots,X_{\tau_h+1},Y)$ and $\alpha^h_{\tau_h+1}(X_1,\dots,X_{\tau_h},J_0^hY)$ are of the same length and perpendicular. Hence, the last curvature ellipse of the polar surface $h$ is a circle.

Finally, observe that the isometric immersion $f$ is written as the composition $f=\Psi_h\circ F,$ where $F\colon U\to UN_{\tau_h}^h$ is the local isometry given by $F(x)=(\pi(x),f(x)), \, x\in U,
$ and $U$ is the saturation of the cross section $L^2\subset M^3$. \qed

\begin{remark}\po{\em
It follows from the computation of the mean curvature of the submanifold $\Psi_h$ in the proof of Theorem \ref{main0}, that the mean curvature is constant by properly choosing the elliptic surface $h.$ Ejiri \cite{Ejiri} proved that tubes in the direction of the second normal bundle of a pseudoholomorphic curve in the nearly K\"ahler sphere $\Sf^6$ have constant mean curvature. Opposed to our case, the index of relative nullity of these tubes is zero.}
\end{remark}

\noindent\emph{Proof of Theorem \ref{main1}:}
Assume that the isometric immersion $f$ is neither minimal nor locally a cylinder. Proposition \ref{sphelliptic} implies that $f$ is spherical. Thus,  from Theorem \ref{main0} we know that for each point on an open dense subset there exist an elliptic surface $h\colon L^2\to \mathbb{Q}^{p-c+4}_c,$  where $c=0$ if $p$ is even and $c=1$ if $p$ is odd, a neighborhood $U$ and a local isometry $F\colon U\to UN^h_{\tau_h}$ such that $f=\Psi_h\circ F.$ In fact,  the elliptic surface $h$ is a polar to $f.$ Moreover, we know that the  curvature ellipses $\mathcal{E}^h_{\tau_h-2}$ and $\mathcal{E}^h_{\tau_h}$  are circles, while the curvature ellipse $\mathcal{E}^h_{\tau_h-1}$ is nowhere a circle.  

Now consider a bipolar surface  $g$ to $f,$ that is, a polar surface to the elliptic surface $h.$ Then it follows from \eqref{eqp} that the curvature ellipse $\mathcal{E}^g_0$ of $g$  is a circle. This means that the bipolar surface is minimal. Furthermore, its first curvature ellipse is nowhere a circle and the second one is a circle. That the isometric immersion $f$ is locally parametrized by \eqref{start} follows from the fact that $f=\Psi_h\circ F$ and
$N_0^g=N_{\tau_h}^h$. \qed

\subsection{Minimal surfaces}

The following proposition provides a way of constructing minimal surfaces in $\R^6$ that satisfy the properties that are required in part (iii) of Theorem \ref{main1}.

\begin{proposition}\po\label{Ricci}
Let $\hat{g}\colon M^2\to \R^6$ be the minimal surface defined by 
$$
\hat{g}=\cos\varphi g_\theta\oplus\sin\varphi g_{\theta+\pi/2},
$$
where $g_\theta, \theta\in [0,\pi), $ is the associated family of a simply connected minimal surface
$g\colon M^2\to \R^3$ with negative Gaussian curvature, and  $\oplus$ denotes the orthogonal sum with respect to an orthogonal decomposition of $\R^6$. If  $\varphi\neq \pi/4,$ then its first curvature ellipse is nowhere a circle and its second curvature ellipse is a circle.
\end{proposition}

Let $g\colon M\to\R^n$ be an oriented minimal surface. The complexified tangent bundle $TM\otimes \mathbb{C}$ is decomposed into the eigenspaces $T^{\prime}M$ and $T^{\prime \prime}M$ of the complex structure $J$, corresponding to the eigenvalues $i$ and $-i.$  The $r$-th fundamental form $\alpha^g_r$, which takes values in the normal subbundle $N_{r-1}^g$, can be complex linearly extended to $TM\otimes\mathbb{C}$ with values in the complexified vector bundle $N_{r-1}^g\otimes \mathbb{C}$  and then decomposed into its $(p,q)$-components, $p+q=r,$ which are tensor products of $p$ differential 1-forms vanishing on $T^{\prime \prime}M$ and $q$ differential 1-forms vanishing on $T^{\prime}M.$ The minimality of $g$ is equivalent to the vanishing of the $(1,1)$-component of the second fundamental form. Hence, the $(p,q)$-components of $\alpha^g_r$ vanish unless $p=r$ or $p=0$. 

It is known (see \cite[Lem. 3.1]{v}) that the curvature ellipse of order $r-1$ is a circle if and only if the $(r,0)$-component of 
$\alpha^g_r$ is isotropic, that is
$$
\<\alpha^g_r(X,\dots,X),\alpha^g_r(X,\dots,X)\>=0
$$
for any $X\in T^{\prime}M,$ where $\<\cdot,\cdot\>$ denotes the bilinear extension over the complex numbers of the Euclidean metric. 
\vspace{2ex}

\noindent\emph{Proof of Proposition \ref{Ricci}:} 
Choose a local tangent  orthonormal frame $e_1,e_2$ such that the shape operator $A$ of $g$ satisfies $AE=k\bar{E},$ where $E=e_1+ie_2$ and $k$ is a positive smooth function. The associated family satisfies 
$g_{\theta_*}=dg\circ J_\theta,$ where $J_\theta=\cos\theta I+\sin\theta J$ and $I$ is the identity endomorphism of the tangent bundle.
Then we have 
\be\label{dghat}
\gh_*E=e^{-i\t}\left(\cos\v g_*E,-i\sin\v g_*E\right). 
\ee

Using the Gauss formula and the fact that the shape operator $A_\theta$ of $g_{\theta}$ is given by $A_\theta=A\circ J_\theta$, we find that the second fundamental form $\hat{\a}$ of $\gh$ satisfies
\be\label{ff2}
\hat{\a}(E,E)=2ke^{-i\t}\left(\cos\v N, -i\sin\v N\right),
\ee
where $N$ is the Gauss map of $g$. It is obvious that 
$\hat{\a}(E,E)$ is not isotropic if $\v\neq \pi/4$, which implies that the first curvature ellipse of $\gh$ is nowhere a circle. 

Differentiating \eqref{ff2} with respect to $E$ and using the Weingarten formula, we obtain 
$$
\tilde{\n}_E\hat{\a}(E,E)=2e^{-i\t}E(k)\left(\cos\v N,-i\sin\v N\right)-2k^2e^{-i\t}\left(\cos\v g_*\bar{E},-i\sin\v g_*\bar{E}\right),
$$
where $\tilde{\n}$ is the connection of the induced bundle of $\gh$. Since $\gh_*E$ and $\gh_*\bar{E}$ span $N_0^\gh\otimes\C,$ the above equation along with \eqref{dghat} yield
$$
\big(\tilde{\n}_E\hat{\a}(E,E)\big)^{N_0^\gh\otimes\C}=-2k^2e^{-2i\t}\cos2\v \gh_*\bar{E}.
$$
It follows using \eqref{ff2} that $N_1^{\gh}\otimes\C={\rm{span}}_{\C}\{\xi,\eta\},$
where $\xi=(N,0)$ and $\eta=(0,iN)$. Then we find that
$$
\big(\tilde{\n}_E\hat{\a}(E,E)\big)^{N_1^\gh\otimes\C}=2e^{-i\t}E(k)\left(\cos\v N,-i\sin\v N\right).
$$
Using the above and since the $(3,0)$-component of the third fundamental form of $\gh$ is given by
$$
\hat{\a}_3(E,E,E)=\big(\tilde{\n}_E\hat{\a}(E,E)\big)^{\left(N^{\gh}_0\otimes\C \oplus N^{\gh}_1\otimes\C\right)^{\perp}},
$$
we obtain
$$ 
\hat{\a}_3(E,E,E)=k^2e^{-i\t}\sin2\v\left(-\sin\v g_*\bar{E},i\cos\v g_*\bar{E}\right).
$$
Thus the (3,0)-component of the third fundamental form of $\gh$ is isotropic, and consequently the second curvature ellipse is a circle. \qed

\section{Submanifolds with constant mean curvature}
In this section, we provide the proofs of the applications of our main results to submanifolds with constant mean curvature.
\vspace{2ex}

\noindent\emph{Proof of Theorem \ref{main2}:} 
The manifold $M^n$ is the disjoint union of the subsets
$$
M_{n-i}=\left\{x\in M^n:\nu(x)=n-i\right\}, \,\ \,\ i=1,2.
$$

Assume that the subset $M_{n-2}$ is nonempty. Then, using Proposition \ref{nu0}  it follows from Theorem \ref{main} for $n\geq4$, or  Theorem \ref{main1} for $n=3$ and $p=1,$ that the isometric immersion $f$ is locally a cylinder over a surface on $M_{n-2}$. 

Suppose that the interior ${\rm{int}}(M_{n-1})$ of the subset $M_{n-1}$ is nonempty. It follows from
 the Codazzi equation that the relative nullity distribution is parallel in the tangent bundle along ${\rm{int}}(M_{n-1})$. Thus, the tangent bundle splits as an orthogonal sum of two parallel orthogonal distributions of rank one and $n-1$ on ${\rm{int}}(M_{n-1})$. By the De Rham decomposition theorem, ${\rm{int}}(M_{n-1})$  splits locally as a Riemannian product of two manifolds of dimension one and $n-1.$  Then, the Gauss equation yields $c=0$. Since the second fundamental form is adapted to the orthogonal decomposition of the tangent bundle, it follows that $f$ is a cylinder over a curve in $\R^{p+1}$ with constant first Frenet curvature (see \cite[Th. 8.4]{da}).
 
Finally, observe that the open subset $V={\rm{int}}(M_{n-1})\cup M_{n-2}$ is dense on $M^n$. \qed
 \vspace{2ex}

In order to proceed to the proofs of the applications of our main results, we need to recall Florit's estimate of the index of relative nullity for isometric immersions with nonpositive extrinsic  curvature. 
The \emph{extrinsic curvature} of an  \ii $f\colon M^n\to\tilde{M}^{n+p}$ for any point $x\in M^n$ and any plane 
$\sigma\in T_xM$ is given by
$$
K_f(\sigma)=K_M(\sigma)-K_{\tilde{M}}(f_*\sigma),
$$
where $K_M$ and $K_{\tilde{M}}$ are the sectional curvatures of $M^n$ and $\tilde{M}^{n+p},$ respectively. 
Florit \cite{f} proved that the index of relative nullity satisfies $\nu\geq n-2p$ at points where the extrinsic curvature of $f$ is nonpositive.
\vspace{2ex}

\noindent\emph{Proof of Corollary \ref{T5}:} We have that the index of relative nullity of $f$ satisfies $\nu\geq n-2.$ Theorem \ref{main2} implies that $c=0$ and, on an open dense subset,  $f$ splits locally as a cylinder over a surface in 
$\R^3$ of constant mean curvature. By real analyticity, the splitting is global.  If $M^n$ is complete, then the surface is also complete with nonnegative Gaussian curvature. That the surface is a cylinder over a circle follows from \cite{ko}. \qed 
\vspace{2ex}

\noindent\emph{Proof of Corollary \ref{C6}:} 
Assume that the hypersurface is nonrigid. Then, the well-known Beez-Killing Theorem (see \cite{da}) implies that the index of relative nullity satisfies $\nu\geq n-2.$ The result follows from Corollary \ref{main2'}. \qed
\vspace{2ex}

\noindent\emph{Proof of Theorem \ref{T7}:} Suppose that the hypersurface is nonminimal.

At first assume that the extrinsic curvature is nonnegative. If $c=0$, a result of Hartman  \cite{hart} asserts that $f(M^n)=\mathbb{S}_R^{k}\times\R^{n-k},$ where $1\leq k\leq n.$ 
If $c=1,$ then $M^n$ is compact by the Bonnet-Myers theorem. According to \cite[Th. 2]{nomichu}, $f$ is totally umbilical. 

In the case of nonpositive extrinsic curvature, the result follows from Corollary \ref{T5}. \qed
\vspace{2ex}

\noindent\emph{Proof of Theorem \ref{T8}:} According to the aforementioned result due to Florit \cite{f},  we have 
$\nu\geq n-4.$ Clearly the manifold $M^n$ is the disjoint union of the subsets
$$
M_{n-i}=\left\{x\in M^n:\nu(x)=n-i\right\}, \,\ \,\ i=1,2,3.
$$

We distinguish the following cases.
\smallskip

\textit{Case I}: Suppose that  the subset $M_{n-4}$ is nonempty. According to Proposition \ref{nu0}, this subset is open. Using \cite[Th. 1]{fz}, we have that on an open dense subset of  $M_{n-4}$ the immersion $f$ is locally  a product $f=f_1\times f_2$ of two hypersurfaces $f_i\colon M^{n_i}\to \R^{n_i+1},i=1,2,$ of nonpositive sectional curvature. The assumption that $f$ has constant mean curvature implies that both hypersurfaces have constant mean curvature as well. Each hypersurface $f_i,i=1,2,$ has index of relative nullity $n_i-2.$
Then it follows from Corollary \ref{main2'} that the submanifold is locally as in part (iii) of the theorem.
\smallskip

\textit{Case II}: Suppose that the interior of the subset $M_{n-3}$ is nonempty. Due to \cite[Th. 1]{fz1}, on an open dense subset of ${\rm{int}}(M_{n-3})$, $f$ is written locally as a composition $f=h\circ F$, where $h=\gamma\times id_{\R^{n-1}} \colon  \R\times\R^n\to\R^{n+2}$ is cylinder over a unit speed plane curve $\gamma(s)$ with nonvanishing curvature $k(s)$ and $F\colon M^n \to\R^{n+1}$ is a hypersurface. The second fundamental form of $f$ is given by
$$
\a^f(X,Y)=h_*\a^F(X,Y)+\a^h\left(F_*X,F_*Y\right),\;\;X,Y\in TM.
$$
From this we obtain $
k\< F_*T,\partial/\partial s\>^2=0 
$
for any $T \in \Delta_f.$ This implies that the height function $F_a=\<F, a\>$ relative to $a=\partial/\partial s$ is constant along the leaves of 
$\Delta_f.$ Then, the mean curvature vector field of $f$ is given by
$$
n\h_f=nH_Fh_*\xi+k\circ F_a\|\grad F_a\|^2\eta,
$$
where $\xi,\eta $ stand for the Gauss maps of $F$ and $h$, respectively. Using that 
$$
\|\grad F_a\|^2=1-\<\xi, a\>^2,
$$
it follows that the mean curvature of $F$ is given as in part (ii) of the theorem.
\smallskip

\textit{Case III}: Suppose that the subset $M_{n-2}\cup M_{n-1}$ has nonempty interior. Then Theorem \ref{main2} implies that the submanifold is locally as in part (i) of the theorem. \qed
\vspace{2ex}

\noindent\emph{Proof of Theorem \ref{T11}:} It follows from \cite[Th. 5.1]{da} that $ \tilde{c}\geq c$ if $n\geq4$.
We distinguish the following cases.
\smallskip

\textit{Case I}: Suppose that $\tilde{c}>c.$ From \cite[Prop. 9]{dt} or \cite[Lem. 8]{O'}, we have that the second fundamental form splits orthogonally and smoothly as 
$$
\alpha^f(\cdot,\cdot)=\beta (\cdot,\cdot)+\sqrt{\tilde{c}-c}\ \<\cdot,\cdot\>\eta,
$$
where $\eta $ is a unit normal vector field and $\beta$ is a flat bilinear form. Thus, the shape operator $A_\xi,$  associated to a unit normal vector field $\xi$ perpendicular to $\eta,$ has $\mathrm{rank}A_\xi\leq1.$ The mean curvature $H$ of $f$ is given by
$$
H^2=\frac{k^2}{n^2} + \frac{\tilde{c}-c}{n},
$$
where $k=\mathrm{trace}A_\xi.$ Obviously, the function $k$ is constant. If $k=0,$ then $f$ is totally umbilical. 

Assume now that $k\neq0.$ Let $X$ be a unit vector field such that $A_\xi X=kX$. The Codazzi equation
$$
(\n_XA_{\eta})T-(\n_TA_{\eta})X=A_{\n_X^\perp\xi}T-A_{\n_T^\perp\xi}X
$$
implies that 
$$
\n_T^\perp\xi=\n_T^\perp\eta=0  
$$
for any $T\in \ker A_\xi$.
Moreover, from the Codazzi equation
$$
(\n_XA_{\xi})T-(\n_TA_{\xi})X=A_{\n_X^\perp\eta}T-A_{\n_T^\perp\eta}X
$$
it follows that
$$
\n_T X=0 \;\;\text{and}\;\;\<\n_X T, X\>=\<\n^\perp_T\xi, \eta\> 
$$
for any $T\in \ker A_\xi$.
Hence the orthogonal distributions $D^1=\mathrm{span}\{X\}$ and $D^{n-1}=\ker A_\xi$ are parallel. By the De Rham decomposition theorem, the manifold splits locally as a Riemannian product $M_{\tilde{c}}^n=M^1\times M^{n-1}$. Consequently, we have $\tilde{c}=0$ and $c=-1$. Clearly $M^{n-1}$ is flat and the second fundamental form is adapted to this decomposition. Then it follows that $f$ is a composition $f=i\circ F$, where $i\colon \R^{n+1}\to \mathbb H^{n+2}$ is the inclusion as a horosphere and  $F\colon M^n_{\tilde{c}} \to \R^{n+1}$ is the cylinder over a circle (see \cite[Th. 8.4]{da}).
\smallskip

\textit{Case II}: Suppose that $c=\tilde{c}.$ It is known that $\nu\geq n-2$ (see Example 1 and Corollary 1 in \cite{Moo}). Then, the result follows from Theorem \ref{main2}.

If $n=3,$ then Theorem \ref{main1} implies that either $c=0$ and $f(M)$ is an open subset of a cylinder over a flat surface $g\colon M^2\to \R^4$ of constant mean curvature, or $c=1$ and $f$ is parametrized by \eqref{start}. In the latter case, it follows from Proposition \ref{sphelliptic} that $f$ is either totally geodesic or elliptic. However, the ellipticity of $f$ implies that the sectional curvature cannot be equal to one. \qed
\vspace{2ex}

\noindent\emph{Proof of Theorem \ref{T9}:} Assume that $f$ is nonminimal. According to Abe \cite{abe}, the index of relative nullity satisfies $\nu\geq n-2.$ Corollary \ref{main2'} implies that the hypersurface is a cylinder over a surface with constant mean curvature. \qed  
\vspace{2ex}

\noindent\emph{Proof of Theorem \ref{T10}:} 
Using \cite[Cor. 2]{fz2} it follows that $\nu\geq n-4.$ The rest of the proof is omitted  since it is similar to the proof of Theorem \ref{T8}.\qed
\vspace{2ex}

The following example produces submanifolds satisfying the conditions in part (ii) of Theorem \ref{T8} or \ref{T10}.
\begin{example}\po{\em Let $F=g\times id_{\R^{n-2}}\colon U\times \R^{n-2}\to \R^{n+1}$ be a cylinder over a rotational surface $g(x,\theta)=(x\cos \theta, x\cos \theta, \varphi(x)), (x,\theta)\in U,$ where $\varphi(x)$ is a smooth function. Consider a cylinder $h=\gamma \times id_{\R^n}$ in $\R^{n+2}$ over a unit speed plane curve $\gamma$ with curvature $k$. Then the isometric immersion $f=h\circ F$ satisfies the conditions in part (ii) of Theorems \ref{T8} and \ref{T10}, with constant constant curvature $H$ and $a=(1,0,\dots,0)$, if the function $\varphi(x)$ solves the ordinary differential equation
$$
\varphi \varphi^{\prime \prime}-1-\varphi^{\prime \, 2}=\pm\varphi\sqrt{(1+\varphi^{\prime \, 2})\left(n^2H^2(1+\varphi^{\prime \, 2})^2-k^2\right)}.
$$
In particular,  $g$ can be chosen as a Delaunay surface and $\gamma$ as the curve with curvature $k=c_0(1+\varphi^{\prime \, 2})$ for a constant $c_0$ such that $0<|c_0|<n|H|$.
}
\end{example}

\bigskip

\noindent Athina Eleni Kanellopoulou\\
University of Ioannina \\
Department of Mathematics\\
Ioannina--Greece\\
e-mail: alinakanellopoulou@gmail.com

\bigskip

\noindent Theodoros Vlachos\\
University of Ioannina \\
Department of Mathematics\\
Ioannina--Greece\\
e-mail: tvlachos@uoi.gr

\end{document}